\documentclass[11pt]{article}

\usepackage{amsmath}
\usepackage{amssymb}
\usepackage{theorem}
\usepackage{amscd}
\usepackage[arrow,curve,matrix,tips,frame]{xy}

\catcode`\@=11 

\renewcommand{\title}[1]{
     \addvspace{3\baselineskip}  
     \begin{center} \LARGE \bf #1
     \end{center}
     \addvspace{2\baselineskip}}   

\renewcommand{\author}[1]{
     \addvspace{-1\baselineskip}  
     \begin{center} \large \sc #1
     \end{center}
     \addvspace{2\baselineskip}}   

\def\section{%
        \@startsection{section}{1}{\z@}%
        {10mm plus 6mm minus 3mm}{\baselineskip}%
        {\normalfont\normalsize\scshape\centering}%
        }

\def\refsection{%
        \@startsection{section}{1}{\z@}%
        {10mm plus 6mm minus 3mm}{\baselineskip}%
        {\normalfont\small\scshape\centering}%
        }

\def\subsection{%
        \@startsection{subsection}{2}{\z@}%
        {\normalbaselineskip}{.35\baselineskip}%
        {\normalfont\normalsize\bf}%
        }

\renewcommand{\paragraph}[1]{{\par\removelastskip\vskip\baselineskip%
         \indent{\scshape{#1}}{\ifperiod.\else\global\periodtrue\fi}%
         \rm \ignorespaces}}

\let\goth=\mathfrak
\let\calligraphy=\mathcal

\def\CC{{\mathbb C}}

\def\PP{{\mathbb P}}

\def\RR{{\mathbb R}}

\def\ZZ{{\mathbb Z}}

\def\Cc{{\calligraphy C}}
\def\Dd{{\calligraphy D}}
\def\Ee{{\calligraphy E}}

\def\Hh{{\calligraphy H}}
\def\Ii{{\calligraphy I}}
\def\Jj{{\calligraphy J}}

\def\Oo{{\calligraphy O}}

\def\Qq{{\calligraphy Q}}

\def\Uu{{\calligraphy U}}

\def\Xx{{\calligraphy X}}
\def\Yy{{\calligraphy Y}}

\def\mmm{{\goth m}}

\def\nnn{{\goth n}}

\def\Ctilde{{\,\widetilde{\!C}}}
\def\Dtilde{{\,\widetilde{\!D}}}

\def\Ftilde{{\,\widetilde{\!F}}}

\def\Ptilde{{\,\widetilde{\!P}}}

\def\Stilde{{\,\widetilde{\!S}}}

\def\Xtilde{{\,\widetilde{\!X}}}
\def\Ytilde{{\,\widetilde{\!Y}}}
\def\Ztilde{{\,\widetilde{\!Z}}}

\def\Xhat{{\,\widehat{\!X}}}

\def\Bl{\operatorname{Bl}}

\def\Div{\operatorname{Div}}

\def\mult{\operatorname{mult}}

\def\Pic{\operatorname{Pic}}
\def\ProjBf{\mathbf{Proj\,}}

\def\supp{\operatorname{supp}}
\def\Sym{\operatorname{Sym}}

\let\lra=\longrightarrow

\def\ie{{\it i.e.}~}
\def\eg{{\it e.g.}~}
\def\inv{^{-1}}

\let\phi=\varphi

\def\Dquad{\hskip 0.6em plus .02em minus .2em}  
\def\Dpar{\belowdisplayskip=0pt\belowdisplayshortskip=0pt\par}

\def\bigpenalty{\interlinepenalty=\@M}
\def\smallpenalty{\interlinepenalty=100}

\newif\ifperiod \periodtrue 

\def\D@makemargins{%
  \labelsep=0pt
  \itemindent=0pt
  \labelwidth=0pt
}
\def\D@restoremargins{%
  \labelsep=5pt
  \itemindent=0pt
  \leftmargin=5mm  
  \labelwidth=\leftmargin \advance\labelwidth by -\labelsep
}
\gdef\th@DthAndSuchtheo{%
  \D@makemargins%
  \def\@begintheorem##1##2{%
  \item[]\hspace\parindent{\scshape ##1~\rm ##2.}\Dquad         
        \D@restoremargins}%
  \def\@opargbegintheorem##1##2##3{\def\next{##3}%
  \item[]\hspace\parindent{\scshape ##1~\rm ##2\ifx\next\empty
  \else\ {\normalfont(##3)}\fi.}         
        \D@restoremargins}%
    }

\gdef\th@DthAndSuchtheostar{%
  \D@makemargins%
  \def\@begintheorem##1##2{%
  \item[]\hspace\parindent{\scshape ##1.}\Dquad     
        \D@restoremargins}%
  \def\@opargbegintheorem##1##2##3{\def\next{##3}%
  \item[]\hspace\parindent{\scshape ##1\ifx\next\empty
  \else\ ##3\fi.}\Dquad         
        \D@restoremargins}%
    }
\gdef\th@DthAndSuchliketheo{
  \D@makemargins%
  \def\@begintheorem##1##2{%
    \@latex@error{likethm: You must provide an argument in square brackets,
    though it may be empty [] !}%
    }%
  \def\@opargbegintheorem##1##2##3{%
        \def\next{##3}\ifx\next\empty\item[\hspace\parindent]\else
        \item[]\hspace\parindent{\scshape \next.}\Dquad\fi
        \D@restoremargins}%
    }

\theoremstyle{DthAndSuchliketheo}
\theorembodyfont{\bigpenalty\itshape}   
\newtheorem{likethm}{}
\theorembodyfont{\rmfamily}  

\theoremstyle{DthAndSuchtheostar}
\theorembodyfont{\bigpenalty\itshape}
\newtheorem{thm*}{\theoname}
\newtheorem{lem*}{\lemmaname}
\newtheorem{pro*}{\propositionname}
\newtheorem{cor*}{\corollaryname}
\newtheorem{conjecture*}{\conjecturename}
\theorembodyfont{\smallpenalty\rmfamily}
\newtheorem{definition*}{\definitionname}
\newtheorem{notation*}{\notationname}
\theorembodyfont{\smallpenalty\rmfamily}
\newtheorem{exa*}{\examplename}
\newtheorem{rem*}{\remarkname}
\theorembodyfont{\smallpenalty\rmfamily}
\newtheorem{definitions*}{\definitionsname}
\newtheorem{examples*}{\examplesname}
\newtheorem{remarks*}{\remarksname}
\theoremstyle{DthAndSuchtheo}
\theorembodyfont{\bigpenalty\itshape}
\newtheorem{thm}{\theoname}[section]
\newtheorem{lem}[thm]{\lemmaname}
\newtheorem{pro}[thm]{\propositionname}
\newtheorem{cor}[thm]{\corollaryname}

\theorembodyfont{\smallpenalty\rmfamily}
\newtheorem{definition}[thm]{\definitionname}

\theorembodyfont{\smallpenalty\rmfamily}
\newtheorem{exa}[thm]{\examplename}
\newtheorem{rem}[thm]{\remarkname}
\theorembodyfont{\smallpenalty\rmfamily}

\newtheorem{remarks}[thm]{\remarksname}

\theoremstyle{DthAndSuchtheo}               
\theorembodyfont{\itshape}

\newcommand{\proof}[1][]{{\par\removelastskip\vskip.6\baselineskip   
        \indent\def\next{#1}{\itshape\proofname\ifx\next\empty\else
        ~\next\fi\ifperiod.\else\global\periodtrue\fi\Dquad}\clubpenalty=5000 
        \rm \ignorespaces}}

\newcommand{\likeproof}[1][]
        {{\par\removelastskip\vskip.6\baselineskip   
        \indent\def\next{#1}\ifx\next\empty\else
        {\itshape{\ignorespaces #1}\ifperiod.\else\global\periodtrue\fi
        \Dquad}\fi\clubpenalty=5000 \rm \ignorespaces}}

\def\qed{{\ifmmode\hskip 6mm plus 1mm minus 3mm{$\square$}\else
\nobreak\hfil\penalty50\hskip1em\null\nobreak\hfil
{\hfill $\square$\parfillskip=0pt\finalhyphendemerits=0\let\par=\endgraf\par}\fi
 \Dpar\penalty-150\vskip.6\normalbaselineskip}}

\def\theoname{Theorem}
\def\lemmaname{Lemma}
\def\propositionname{Proposition}
\def\notationname{Notation}
\def\corollaryname{Corollary}
\def\conjecturename{Conjecture}
\def\remarkname{Remark}
\def\remarksname{Remarks}
\def\examplename{Example}
\def\examplesname{Examples}
\def\definitionname{Definition}
\def\definitionsname{Definitions}
\def\notationname{Notation}

\def\proofname{Proof}

\catcode`\@=12 

\newcommand{\arhook}[1]{\ar@^{(->}[{#1}] \ar[{#1}]}

\def\Deltatilde{{\,\widetilde{\!\Delta}}}
\def\bEe{\boldsymbol{\Ee}}
\def\bHh{\boldsymbol{\Hh}}
\def\bYy{\boldsymbol{\Yy}}

\def\etilde{{\,\widetilde{\!e}}}

\newcommand{\arrow}[1]{\,\stackrel{#1}{\longrightarrow}\,}

\def\Bs{\operatorname{Bs}}
\def\elm{\operatorname{elm}}
\def\Ch{\mathcal{C}h}
\def\sff{\mbox{\sf f}}

\newcounter{icounter}
\newenvironment{ilist}
  {\begin{list}{\upshape (\roman{icounter})}
    {\setlength{\leftmargin}{25pt}
     \setlength{\labelwidth}{20pt}
     \setlength{\labelsep}{5pt}
     \setlength{\topsep}{1pt plus 1pt minus 1pt}
     \setlength{\parsep}{0pt}  
     \setlength{\itemsep}{0pt plus .5pt minus 1pt}
     \usecounter{icounter}} }
  {\end{list}}

\begin{document}
\title{Rationality properties of manifolds containing quasi-lines}
\author{Paltin Ionescu and Daniel Naie}

\begin{flushleft}
\end{flushleft}

\begin{abstract}
Let $X$ be a complex, rationally connected, projective manifold.
We show that $X$ admits a modification $\Xtilde$ that contains a
quasi-line, \ie a smooth
rational curve whose normal bundle is a direct 
sum of copies of $\Oo_{\PP^1}(1)$. 
For manifolds containing quasi-lines, a sufficient condition of
rationality is exploited: There is a unique quasi-line from a given
family passing through two general points. 
We define a numerical birational invariant, $e(X)$, 
and prove that $X$ is rational if and only if $e(X)=1$.
If $X$ is rational, there is a modification $\Xtilde$ which is
strongly-rational, \ie contains an open subset isomorphic to an open 
subset of the projective space whose complement is at least 
$2$-codimensional. We prove that strongly-rational 
varieties are stable under smooth, small deformations. 
The argument is based on a convenient caracterization of 
these varieties.

Finally, we relate the previous results and formal geometry. This
relies on $\etilde(X,Y)$, a numerical invariant of a given quasi-line
$Y$ that depends only on the formal completion $\Xhat|_Y$. As
applications we show various instances in which $X$ is determined by
$\Xhat|_Y$. We also formulate a basic question about the birational
invariance of $\etilde(X,Y)$.
\end{abstract}

\section*{Introduction}
 
Classical examples of rational projective manifolds are given by
usually elementary, sometimes ingenious, geometric constructions of
linear systems, yielding birational maps (\eg projections from
subvarieties). Related to the L\" uroth problem in dimension at least 
three, several fairly sophisticated techniques for proving non-rationality
of some Fano manifolds have been developped (see \eg \cite{IP}). Using
deformation theory of rational curves, Koll\'ar, Miyaoka and Mori
introduced in \cite{KMM} the very useful class of rationally connected
varieties, generalizing the classes of both rational and Fano
manifolds. 
Rational connectedness admits several convenient characterizations and is
invariant under deformations and birational isomorphism. It is therefore
natural to try to understand rationality within the larger class of
rationally connected manifolds. 
 
Let $X$ be a complex projective manifold of dimension $n \geq 2$. 
$X$ is {\it rationally connected} if two general points of
it may be joined by a rational curve. Equivalently, $X$ contains a smooth 
rational curve with ample normal bundle, see \cite{KMM}, \cite{Ko}.
A smooth rational curve $Y \subset X$ is called a {\it quasi-line} 
(see \cite{BBI}) if its normal bundle is 
isomorphic to $\bigoplus_1^{n-1} \Oo_{\PP^1}(1)$. 
$X$ is called {\it strongly-rational}
(see \cite{BBI}) if there exists a birational map 
$\varphi : X \dashrightarrow \PP^n$ 
which is an isomorphism from an open subset $U$ onto an open subset $V$, 
whose complement in $\PP^n$ is at least $2$-codimensional. 
Note that strongly-rational manifolds contain quasi-lines (the pull-back 
of a line contained in $V$ gives rise to a quasi-line on $X$).
Therefore we have the diagram:
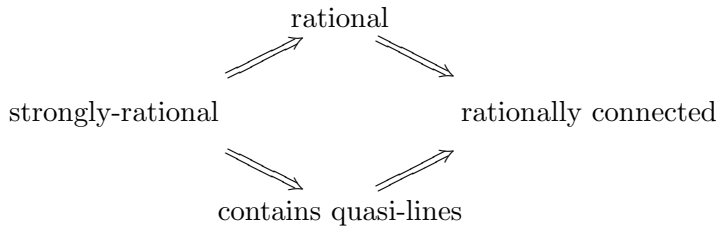
\[
\begin{xy} <0.5cm,0cm>:
 (-3,0)*+!R{\shortstack{strongly-rational}},
 (0,2.1)*+!D{\shortstack{rational}},
 (3,0)*+!L{\shortstack{rationally connected}}, 
 (0,-2.1)*+!U{\shortstack{contains quasi-lines}},
    \ar@{=>} (-3,1);(-1,2)
    \ar@{=>} (1,2);(3,1)
    \ar@{=>} (-3,-1);(-1,-2)
    \ar@{=>} (1,-2);(3,-1)
\end{xy}
\]

In the first section we quote from \cite{I} and \cite{Io} two rather
general classes of examples of rational manifolds. Moreover, we prove a
new rationality criterion, Theorem \ref{1.3}.
 
 In Section \ref{s2} we show that any rationally connected manifold, after
being suitably blown-up, contains quasi-lines, Theorem \ref{2.3}. 
Note that, in a similar vein, by \cite{Hi}, rational manifolds become 
strongly-rational after suitable blowing-ups. The proof of Theorem
\ref{2.3} applies to show the existence of almost-lines, \ie
quasi-lines $Y$ such that $D\cdot Y=1$ for some divisor $D$ on
$X$. This completes Theorem 2.1 from \cite{BBI}. 
 
In Section \ref{s3} we use quasi-lines to characterise rationality and
to define, for each rationally connected manifold $X$, a birational
numerical invariant, denoted $e(X)$. We first introduce
and compute for some examples, the number $e(X,Y)$ of quasi-lines from
a given family that pass through two general points of $X$; for
instance, when $X$ is a smooth cubic threefold in $\PP^4$ and $Y$ is a
general conic, $e(X,Y)=6$, see Proposition \ref{e0eForConic}. Then,
$e(X)$ represents the minimum among $e(X',Y')$, where $X'$ is obtained
from $X$ by a sequence of blowing-ups with smooth centers. In Theorem
\ref{birEx}, we prove that $X$ is rational if and only if
$e(X)=1$. However, note that $e(X)$ seems to be very difficult to compute.
In order to get the  rationality via quasi-lines, the key
is to show that $e(X,Y)=1$ for a certain quasi-line $Y$, see
Proposition \ref{e0e}. 

Section \ref{s4} contains a convenient characterization of
strongly-rational manifolds, Theorem \ref{4.2}. As a consequence we
show in Theorem \ref{4.5} that strongly-rational manifolds are stable
with respect to small deformations. Note that such an
invariance property is not expected to hold for rational manifolds. 
 
Section \ref{formalGeometry} relates the preceding results to formal
geometry. To each quasi-line $Y\subset X$ we associate a ``local''
invariant denoted $\etilde(X,Y)$. It depends only on the formal
completion $\Xhat|_Y$. Theorem \ref{eEtildeB} shows that
$e(X,Y)=\etilde(X,Y)\cdot b(X,Y)$, where
$b(X,Y)=[K(\Xhat|_Y):K(X)]$. Here $K(X)$ is the field of rational
functions on $X$ and $K(\Xhat|_Y)$ is the field of formal rational
functions of $X$ along $Y$, see \cite{HM}. As applications of this
formula we give various examples of instances when the formal
completion of $X$ along the quasi-line $Y$ determines $(X,Y)$ up to
isomorphism. They include the case mentioned above of a general conic
on a cubic threefold, Corollary \ref{conicOnCubic}.

In the last section we address the basic question: Does the local
invariant $\etilde(X,Y)$ depend only on the field $K(\Xhat|_Y)$?
A positive answer would have nice consequences, \eg a completely new
proof of the non-rationality of the smooth cubic threefold in $\PP^4$.

In the Appendix, we show via a toric calculation, that a certain useful
property of quasi-lines does not hold in general.

We shall work over the field of complex numbers. Unless otherwise
stated, we follow the usual conventions and notation in Algebraic
Geometry (see \eg \cite{Ha}).

\section{Some rational varieties}

 Let $X\subset \PP ^N$ be a projective manifold of dimension $n \geq 2$. In
dimension two, the famous Castelnuovo criterion characterizes rationality
by the vanishing of two numbers which are birational invariants of $X$;
in particular, rationality and rational connectedness are equivalent.
For $n \geq 3$, deciding the rationality of $X$ may be a quite difficult
problem. Rationally connected manifolds, which are easier to understand,
form a much larger class than rational ones.
 
 Many examples of rational manifolds come from more precise 
biregular classification statements. We would like to exemplify this
principle by two rather general results. To state them we recall some
numerical invariants of $X$.
 
We denote by $g$ the {\it sectional genus} of $X$, that is the genus 
of the curve got by intersecting $X$ with $n-1$ general hyperplanes in 
$\PP ^N$. We let $d$ be the {\it degree} of $X$. Finally, 
let $q=:h^1(X, \Oo_X)$ be the {\it irregularity} of $X$.
 
\begin{thm} \label{1.1}
Assume that $d \geq 2g-1$ and $ q=0$. Then $X$ is rational, unless $X$ is
a cubic hypersurface, $n\geq 3$.
\end{thm}

This statement is a consequence of the precise biregular classification,
given in \cite{I}, Corollaries 8, 9 and 10, of all manifolds satisfying
conditions $d \geq 2g-1$, $n \geq 3$. The classification, due to Fujita,
of the so called ``del Pezzo manifolds", is also used. It corresponds to the
case $g=1$, which includes the exception in the statement of Theorem 1.1.
Note that the bound $d \geq 2g-1$ is sharp. Indeed, the quartic threefold
in $\PP ^4$ and the complete intersection of a quadric and a cubic in
$\PP ^5$ satisfy $d=2g-2$. However, they are known, as well as the cubic
threefold of $\PP ^4$, to be non-rational. See \eg \cite{IP} for a 
discussion of these very delicate results.
 
The previous theorem shows that, for fixed sectional genus, regular manifolds
of ``high" degree are rational. On the other hand, many examples of manifolds 
of ``small" degree are known to be rational. Moreover, note that deciding
the rationality property is particularly difficult when $X$ is a Fano manifold
with $b_2(X)=1$. In this direction, we quote the following recent result:

\begin{thm} \label{1.2} {\rm (cf.~\cite{Io})}
Let $X\subset\PP^N$ be non-degenerate and assume that $d\leq N$. 
Then one of the following holds:
\begin{ilist}
  \item $X$ is Fano and $b_2(X)=1$, or
  \item $X$ is rational.
\end{ilist}
\end{thm}

The bound $d \leq N$ is clearly the best possible one. A hypersurface 
of degree $N+1$ is neither rational, nor Fano. Again, the rationality comes
a posteriori, using a classification result. In fact, manifolds as in (ii)
may be completely described: There are 6 infinite series and 14 ``sporadic"
examples, see \cite{Io}.

 Next we prove a result that allows one to deduce the 
rationality directly from the existence of a suitable rational submanifold
of $X$. To the best of our knowledge, this theorem seems to have 
been overlooked in the classical literature. Our proof depends on 
Hironaka's desingularisation theory from \cite{Hi} and on basic 
properties of rationally connected manifolds (cf.~\cite{KMM}).

\begin{thm} \label{1.3}
Let $X$ be a projective variety and $|D|$ a complete linear system of 
Cartier divisors on it. Let $D_1,\ldots,D_s \in |D|$ and put
$W_i=:D_1\cap\cdots\cap D_i$ for $1 \leq i \leq s$. Assume that 
for all $i$, $W_i$ is smooth, irreducible and has dimension $n-i$. 
Assume moreover that there is a divisor $E$ on $W=:W_s$ 
and a linear system $\Lambda \subset |E|$ such that
\begin{ilist}
  \item $\varphi_\Lambda:W \dashrightarrow \PP ^{n-s}$ is birational, and 
  \item $|D|_W - E| \neq \emptyset$.
\end{ilist}
Then $X$ is rational.
\end{thm}
\proof
We proceed by induction on $s$. We explain the case $s=1$, the general 
induction step being completely similar. So, let $W \in |D|$ be a smooth,
irreducible Cartier divisor such that 
$\varphi_\Lambda:W \dashrightarrow \PP ^{n-1}$
is birational for $\Lambda \subset |E|$, $E\in \Div(W)$ and $|D|_W - E|\neq
\emptyset$. Note that $W$ is contained in the smooth locus of $X$. So,
replacing $X$ by its desingularisation, we may assume $X$ to be smooth.
As $W$ is rational, it is in particular rationally connected; so by \cite{KMM},
there is some smooth rational curve $Y \subset W$ with ample normal bundle.
We have $Y \cdot E > 0$ since E moves and $Y \cdot(D|_W - E) \geq 0$ 
by condition (ii).
It follows that $Y \cdot D > 0$. Looking at the standard exact sequence of
normal bundles, we get that $N_{Y|X}$ is ample. So, again by \cite{KMM},
$X$ is rationally connected and, in particular, $q(X)=h^1(X,\Oo_X)=0$.
The standard exact sequence:
\[
  0 \lra \Oo_X 
    \lra \Oo_X(D) 
    \lra \Oo_{W}(D) 
    \lra 0,
\]
shows that $\dim |D|=\dim |D|_W| + 1 \geq \dim |E| + 1 \geq n$. 
Choose a pencil
$(W,W') \subset |D|$, containing $W$, such that $W'|_W=E_0 + E_1$, 
with $E_1\in \Lambda$ and $E_0 \geq 0$. Now, by the theory in
\cite{Hi}, we may use blowing-ups with smooth centers contained in 
$W\cap W'$, such that, after taking the proper transforms  
of the elements of our pencil, to get:

 (a) $\supp (E_0)$ has normal crossing;

 (b) $\Lambda$ is base-points free (so $\varphi:W \to \PP ^{n-1}$ is a
birational morphism).\\
Next, by blowing-up the components of $\supp(E_0)$, we may also suppose that
$E_0 = 0$, \ie $D|_W$ is linearly equivalent to $E$. Now, using the
previous standard exact sequence and the fact that $q(X) = 0$, it follows that
$|D|$ is base-points free. But we have that $D^n=(D|_W)^{n-1}_W=1$, 
so $\varphi_{|D|}$ is a birational morphism onto $\PP^n$. 
\qed

\begin{exa} 
Let $X \subset \PP^{n+d-2}$ be a non-degenerate projective variety of
dimension $n \geq 2$ and degree $d \geq 3$. Then $X$ is rational, unless
it is a smooth cubic hypersurface, $n\geq 3$.

Indeed, 
we may assume $X$ to be smooth; otherwise, use a projection from a singular
point. We may also suppose that $X$ is linearly normal (if not, use again
a projection from one of its points). One sees easily that such a linearly
normal, non-degenerate manifold $X \subset \PP^{n+d-2}$ has anticanonical 
divisor linearly equivalent to $n-1$ times the hyperplane section, \ie
they are exactly the so called ``classical del Pezzo manifolds''. They were
classified by Fujita in a series of papers, see \eg \cite{IP} for a survey
of his argument. As Fujita's proof is quite long and difficult, we show
how Theorem 1.3 above may be used to prove directly the rationality of
$X$ if $d \geq 4$.
 Consider the surface $W$ got by intersecting $X$ with $n-2$ general
hyperplanes. Note that $W$ is a non-degenerate, linearly normal surface of
degree $d$ in $\PP ^d$, so it is a del Pezzo surface. As such, $W$ is known
to admit a representation $\varphi:W \to \PP ^2$ as the blowing-up of
$9-d$ points (in general position). Let $L\subset W$ be the pull-back via
$\varphi$ of a general line in $\PP ^2$. It is easy to see that $L$ is a 
cubic rational curve in the embedding of $W$ into $\PP ^d$. So, for $d \geq 4$, 
$L$ is contained in a hyperplane of $\PP ^d$. This shows that the conditions of
Theorem 1.3 are fulfilled for $X$, $|D|$ being the system of
hyperplane sections. We also see that Theorem 1.3 is sharp, as the 
previous argument fails exactly for the case of cubics.
\end{exa}

In the remaining of this section we slightly generalize the fibration
Theorem 1.12 in \cite{IoVo}. As a consequence we get a rationality
criterion, Corollary \ref{c2}; it was this criterion that led us to
formulate Theorem \ref{1.3}. It will be  convenient to refer to a
couple $(X,Y)$, where $Y$ is a smooth rational curve with ample normal
bundle, as to a model, cf.~\cite{IoVo}. Talking about a model $(X,Y)$,
we shall often replace $(X,Y)$ by $(X,Y')$, where $Y'$ is a
deformation of $Y$ (we write $Y' \sim Y$). A morphism of models,
$(X,Y) \to (X',Y')$, is a morphism $X \to X'$ that maps $Y$
isomorphically to $Y'$.

Firstly, we recall the above mentioned fibration theorem:

\begin{thm}[cf.~\cite{IoVo}, 1.12] \label{thIoVo}
Let $(X,Y)$ be a model  and let $D$ be a divisor such that 
$D \cdot Y = 1$ and $\dim |D|=:s \geq 1$. Then there 
exists $Y' \sim Y$ and a diagram of models
\[
\UseTips
 \newdir{ >}{!/-5pt/\dir{>}}
\xymatrix{   
  (Z,\Ytilde') \arhook{r} & 
    (\Xtilde,\Ytilde') \ar[r]^\varphi \ar[d]_\sigma & (\PP^s,l)  \\
  & (X,Y')}
\]
such that
\begin{ilist}
  \item $\sigma$ is a sequence of blowing-ups,
  \item $\varphi$ is surjective, with connected fibres,
  \item any smooth fibre of $\varphi$ is rationally connected,
  \item $l \subset \PP^s$ is a line and $Z = \varphi^{-1}(l)$ is smooth, and
  \item $\Ytilde'$ is a section for $\varphi |_Z$.
\end{ilist}
\end{thm}

Next, we generalize it to the case when 
$D\cdot Y\geq 2$. To see this, we observe the behaviour of the normal
bundle of the curve when $X$ is blown-up at a point lying on the
curve. 

\begin{lem} 
Let $C \subset X$ be a smooth curve, $p \in C$ a point 
and $\sigma : \Xtilde \to X$ the blowing-up of $X$ at $p$.
If $\Ctilde$ is the strict transform of $C$, then 
\[
N_{\Ctilde| \Xtilde} \simeq 
  \sigma^\ast (N_{C|X} \otimes \Oo_{C}(-p)).
\]
\end{lem}
\proof
Let $\{U_\alpha\}$ be a covering of $C$ with open subsets of $X$, 
$p \in U_0$. Let $(u_1^\alpha, \ldots ,u_n^\alpha)$ 
be local coordinates on $U_\alpha$ such that 
$u_1^\alpha, \ldots, u_{n-1}^\alpha$ are local equations 
for $C$, and $u_n$ is a local equation for $p$ along $C$. If 
\[
  u_i^\alpha = \sum_{j=1}^{n-1} h_{ij}^{\alpha \beta} u_j^\beta 
    \qquad \text{for every} \quad i = 1, \ldots , n-1
\]
on $U_\alpha \cap U_\beta$, then 
$c^{\alpha \beta} = (h_{ij}^{\alpha \beta}|_C)$ represents 
the transition function for $N_{C|X}$, from $U_\alpha$ to 
$U_\beta$. The open covering $\{U_\alpha\}$ induces an open covering 
$\{V_\alpha\}$ of $\Ctilde$: for $\alpha \neq 0$, 
$V_\alpha = U_\alpha$. For $\alpha = 0$, the open subset $V_0$ 
is an open subset with local coordinates 
$(\xi_1, \ldots, \xi_{n-1}, u_n^0)$ such that $u_i^0 = u_n^0 \xi_i$. 
It follows that $N_{\Ctilde|\Xtilde}$
is given by the transition functions $c^{\alpha \beta}$, if 
$\alpha \neq 0$ and $\beta \neq 0$, and by 
$\frac{1}{u_n^0}c^{0\beta}$ if not. 
But the constant function $1$ and $\frac{1}{u_n^0}$ are transition
functions of $\Oo_C(-p)$ relative to the covering $\{U_\alpha\}$, 
hence the result.
\qed

Now we can rephrase Theorem \ref{thIoVo} as follows:

\begin{thm} 
Let $(X,Y)$ be a model such that $N_{Y|X}=\bigoplus_1^{n-1}\Oo_Y(a_j)$
with $a_1 \leq \cdots \leq a_{n-1}$, and let $D$ be a divisor such that 
$D \cdot Y =: d > 0$ with $a_1 \geq d$ and $\dim |D| \geq d$. Then, 
there is $\Xtilde$ a blow-up of $X$ and a diagram of models
\[
  (Z,\Ytilde') \, \, \hookrightarrow \, \,
  (\Xtilde,\Ytilde')\arrow{\varphi}
  (\PP^{\dim|D|-d+1},l)
\]
such that
\begin{ilist}
  \item $\varphi$ is surjective, with connected fibres,
  \item any smooth fibre of $\varphi$ is rationally connected,
  \item $Z = \varphi^{-1}(l)$ is smooth, and
  \item $\Ytilde'$ is a section for $\varphi|_Z$.
\end{ilist}
\end{thm}
\proof
We may suppose that $|D|$ is free from fixed components
and that $Y$ does not meet the base locus 
of $|D|$. Indeed, for the latter, 
a general deformation of $Y$ avoids a closed subset of codimension
$\geq 2$ and in the decomposition of
its normal bundle, $a_1 \geq d$ (the function $-\min a_j$ is
upper-semicontinuous, see \cite{Ko}, Lemma II.3.9.2).
We continue by blowing up $d-1$ points on $Y$. We take $X'$ to be 
the new variety and $Y'$ the strict transform of $Y$; 
then by the above lemma $(X',Y')$ is a model. 
Moreover, the divisors linearly equivalent to $D$ through 
the $d-1$ points determine a linear system  $|D'|$ on $X'$. 
We have $D'\cdot Y'=1$ and $\dim |D'| \geq 1$, so Theorem \ref{thIoVo}
applies.
\qed

\begin{cor} \label{c1}
Let $(X,Y)$ be a model with 
$N_{Y|X}= \bigoplus_1^{n-1} \Oo_Y(a_j)$, $a_1 \leq\cdots\leq a_{n-1}$ 
and let $D$ be a divisor such that $0<D \cdot Y=:d \leq a_1$ and 
$\dim |D|\geq n+d-1$. Then $X$ is rational.
\end{cor}

Finally, we state the following corollary from which Theorem \ref{1.3}
stemmed. 

\begin{cor} \label{c2}
Let $X$ be a projective variety and $|D|$ a complete linear system of 
Cartier divisors on it. Let $D_1,\ldots,D_s \in |D|$ and put
$W_i=:D_1\cap\cdots\cap D_i$ for $1 \leq i \leq s$. Assume that 
for all $i$, $W_i$ is smooth, irreducible and has dimension $n-i$. 
Assume moreover that there is a divisor $E$ on $W=:W_s$ 
and a smooth rational curve $Y\subset W$ such that:
\begin{ilist}
  \item $N_{Y|W}\simeq\bigoplus_1^{s-1}\Oo_{\PP^1}(a_i)$,
  $0<d\leq a_1\leq\cdots\leq a_{s-1}$, where $d=:E\cdot Y$, 
  \item $\dim |E|\geq s+d-1$, and
  \item $|D|_W - E| \neq \emptyset$.
\end{ilist}
Then $X$ is rational.
\end{cor}
\proof
We proceed as in the proof of Theorem \ref{1.3}. After suitable
blowing-ups we get $D\cdot Y=d$ and $\dim|D|\geq n+d-1$. So the above
corollary applies.
\qed

\section{Existence of quasi-lines and a first application} \label{s2}

 In this section we show that a rationally connected 
manifold, up to blowing it up along smooth subvarieties, 
contains quasi-lines. The proof depends on the following considerations
about elementary transforms, which may be of some independent interest.
 
 If $M$ is a smooth variety and $V \to M$ a vector bundle,
then we shall use the classical convention for the 
projective space that is most suitable for the present work.
Accordingly, the associated projective bundle of $V$ will be 
$P(V) = \ProjBf (\Sym V^\ast)$.
 
 Let $C$ be a smooth curve, $V \to C$ a vector bundle of rank $n$ 
and $\pi : P(V) \to C$ the corresponding projection. 
If $F$ is a fibre of $P(V)$, $\pi(F) = c$,
and $L \subset F$ a hyperplane, the elementary transform of $P(V)$ 
with center $L$, denoted by $\elm_L P(V)$, is the projective bundle 
$P'$ over $C$ constructed as follows: 

1) Denote by $\Ptilde$ the blow-up of $P=P(V)$ along 
$L$, and by $\sigma$ the projection from $\Ptilde$ to $P$.
The exceptional divisor $E$ of $\Ptilde$ is $P(N_{L|P})$, 
a $\PP^1$-bundle over $L$. 
The fibre $(\pi \circ \sigma)^{-1}(c)$ is the sum of two effective
Cartier divisors
\[
(\pi \circ \sigma)^{-1}(c) = \Ftilde + E.
\]
They intersect in the hyperplane of $\Ftilde$ that corresponds 
to $L$. On $E$, this intersection is the exceptional divisor, when 
$E$ is seen as the blow-up of $\PP^{n-1}$ at a point.

2) The normal bundle of $\Ftilde$ in $\Ptilde$ is 
$\Oo_{\Ftilde}(-1)$, hence there is a contraction 
$\sigma' : \Ptilde \to P'$ that sends $\Ftilde$ to a point: 
\[
\UseTips
 \newdir{ >}{!/-5pt/\dir{>}}
 \xymatrix{
               & \Ptilde \ar[ld]_\sigma \ar[rd]^{\sigma'} & \\
 P \ar[rd]_\pi &                           & P' \ar[ld]^{\pi'} \\ 
               & C                         &}
\]
$P'$ maps to $C$ with all fibers isomorphic to $\PP^{n-1}$. 
It follows that $P'$ is a projective bundle.

\bigskip

The construction of $\elm_L P(V)$ is a generalization of the 
elementary transforms of geometrically ruled surfaces. 
In this case, if the base curve is the projective line, 
the elementary transform can be described more precisely: The 
Hirzebruch surface $P(\Oo_{\PP^1} \oplus \Oo_{\PP^1}(d))$ is 
transformed either in the surface 
$P(\Oo_{\PP^1} \oplus \Oo_{\PP^1}(d+1))$, or the surface 
$P(\Oo_{\PP^1} \oplus \Oo_{\PP^1}(d-1))$, depending on whether 
or not the point $L$ lies on the distinguished section (see \cite{Ha}). 
The next proposition is an analogous result in 
arbitrary dimensions.

\begin{pro} \label{2.1}
Let $V = \Oo_{\PP^1}(a_1) \oplus \cdots \oplus \Oo_{\PP^1}(a_n)$, 
with $a_1 \leq \cdots \leq a_n$. If $L$ is a general hyperplane 
in a fibre of $P(V) \to \PP^1$, then $\elm_L P(V) = P(V')$, where 
\[
  V' = \Oo_{\PP^1}(a_1) \oplus \cdots 
       \oplus \Oo_{\PP^1}(a_{n-1}) \oplus \Oo_{\PP^1}(a_n-1).
\]
\end{pro}
\proof
If $W$ is a sub-bundle of $V$, with quotient bundle $Q$, there is 
a canonical embedding $i$ of $P(W)$ in $P(V)$. For $W$ of rank $n-1$, 
$P(W)$ is an effective divisor. Taking $\pi_\ast$ on the exact
sequence 
\[
  0 \lra \Oo_V(1) \otimes \Oo_V(-P(W)) 
    \lra \Oo_V(1) 
    \lra \Oo_V(1) \otimes \Oo_W 
    \lra 0,
\]
we have 
\[
  0 \lra \pi_\ast(\Oo_V(1) \otimes \Oo_V(-P(W)) 
    \lra V^\ast 
    \lra W^\ast 
    \lra 0.
\]
It follows that $\Oo_V(P(W)) \otimes \Oo_V(-1) \simeq \pi^\ast Q$, 
and restricting to $P(W)$, that 
\[
  \Oo_W(P(W)) \otimes \Oo_W(-1) \simeq (\pi \circ i)^\ast Q.
\]

We shall call splitting divisors, 
the divisors $\Delta_i$, $1 \leq i \leq n$, corresponding to 
the sub-bundles 
\[
  W_i = \Oo_{\PP^1}(a_1) \oplus \cdots 
        \oplus \Oo_{\PP^1}(a_{i-1}) \oplus \Oo_{\PP^1}(a_{i+1}) 
        \oplus \cdots \oplus \Oo_{\PP^1}(a_n) \subset V.
\]
These $n$ divisors have empty set-theoretic intersection and 
\begin{equation} \label{SelfIntersections}
  (\Delta_i)^n = c_1^n(\Oo_V(1) \otimes \Oo_V(a_iF))
               = - (a_1 + \cdots + a_n) + na_i,
\end{equation}
for any $1 \leq i \leq n$.
It is obvious that giving a projective bundle over $\PP^1$ is
equivalent to giving $n$ $1$-codimensional projective sub-bundles 
with an empty intersection. 
The splitting type of the bundle can be restored, up to tensoring 
with a line bundle, from the self intersection numbers 
(\ref{SelfIntersections}). 

\medskip

\noindent 
{\it Claim:} 
Any sufficiently general hyperplane $L \subset F$ is cut out by a 
projective sub-bundle $\Delta$ linearly equivalent to $\Delta_n$.

It is sufficient to show that the restriction 
$H^0(P(V), \Oo_V(\Delta_n)) \to H^0(F, \Oo_F(1))$ is surjective. 
The dimension
\[
\begin{split}
  h^1(P(V), \Oo_V(\Delta_n-F)) 
    &= h^1(P(V), \Oo_V(1) \otimes \Oo_V((a_n-1)F))             \\
    &= h^1(\PP^1, \pi_\ast (\Oo_V(1) \otimes \Oo_V((a_n-1)F))) \\
    &= h^1(\PP^1, \makebox{$\bigoplus$}_{i=1}^n \Oo_{\PP^1}(-a_i+a_n-1)F))\\
    &= h^0(\PP^1, \makebox{$\bigoplus$}_{i=1}^n \Oo_{\PP^1}(a_i-a_n-1)F))
\end{split}
\]
vanishes, since $a_i \leq a_n$ for all $i$, and the surjection follows. 
If $p_n \in F$ is the intersection of the fibre $F$ with the 
section that corresponds to $\Oo_{\PP^1}(a_n) \hookrightarrow V$,
then $\Delta_n$ does not pass through $p_n$, contrary to the other 
$\Delta_i$'s.

\medskip

Let $L$ be a general hyperplane in $F$, with $p_n \notin L$, 
let $\Delta \sim \Delta_n$ be the projective sub-bundle that 
corresponds to $L$, and let $P' = \elm_L P$.  
We reconstruct the vector bundle corresponding to $P'$ from the 
$n$ $1$-codimensional sub-bundles 
$\Delta'$, $\Delta'_1, \ldots, \Delta'_{n-1}$, where 
$\Delta' = \sigma'(\Deltatilde)$, with $\Deltatilde$ the 
strict transform of $\Delta$ on $\Ptilde$, and the same for  
$\Delta'_i$, $i=1, \ldots, n-1$. 
Let $b_1, b_2, \ldots, b_n$ be $n$ integers such that 
$(\Delta'_i)^n = -(b_1 + \cdots + b_n) + nb_i$
and $(\Delta')^n = -(b_1 + \cdots + b_n) + nb_n$.
Then 
$P' = P(\Oo_{\PP^1}(b_1) \oplus \cdots \oplus \Oo_{\PP^1}(b_n))$, 
and since for every $1 \leq i \leq n-1$, 
$\sigma^\ast \Delta_i = \Deltatilde_i$ and 
$(\sigma')^\ast \Delta'_i = \Deltatilde_i + \Ftilde$, 
we have that
\[
  (\Delta_i)^n = (\Deltatilde_i)^n 
               = ((\sigma')^\ast \Delta'_i - \Ftilde)^n 
               = (\Delta'_i)^n + (-1)^n \Ftilde^n 
               = (\Delta'_i)^n - 1.
\]
These relations together with (\ref{SelfIntersections})  
provide a linear system of $n-1$ equations, 
$n$ unknowns and of rank $n-1$. Now $b_i = a_i$ for $1 \leq i \leq n-1$, 
and $b_n = a_n-1$ is one solution, the others being obtained from this
one by translations.
\qed

\begin{lem} \label{2.2}
Let $X$ be a smooth projective $n$-fold and $Y \subset X$ be a smooth 
rational curve with normal bundle 
$N_{Y|X} = \Oo_Y(a_1) \oplus \cdots \oplus \Oo_Y(a_{n-1})$, 
where $a_1 \leq \cdots \leq a_{n-1}$. If $Z \subset X$ is a 
general, smooth, $2$-codimensional subvariety intersecting 
$Y$ in a point, $X' = \Bl_Z X$ and 
$\Ytilde$ the strict transform of $Y$, then 
\[
  N_{\Ytilde|X'} = 
    \Oo_{\Ytilde}(a_1) \oplus \cdots \oplus \Oo_{\Ytilde}(a_{n-1}-1).
\]
\end{lem}
\proof
We need to compare the normal bundles of $Y$ and
$\Ytilde$. Accordingly, we first look for a comparison
of the corresponding projective bundles. These are exceptional divisors 
of the blow-ups of $X$ and $X'$ along $Y$ and $\Ytilde$, 
respectively. Throughout the proof, the exceptional divisor of 
the blow-up along the subvariety $S$ will be denoted by 
$P_S$, and the strict transform of a subvariety $S$ on a blow-up 
by $\Stilde$. Hence, in the next diagram, we want to 
compare the exceptional divisors $P_Y\subset\Bl_Y X$ and 
$P_{\Ytilde}\subset\Bl_{\Ytilde} X'$.
\[
\begin{CD}
  P_{\Ytilde} \subset \Bl_{\Ytilde} X' @>{\sigma'}>> X' \\
    &&                                       @VV{\rho}V    \\
  P_Y \subset \Bl_Y X                        @>\sigma>> X
\end{CD}
\]

\medskip

\noindent
{\it Claim:} $P_{\Ytilde} = \elm_L P_Y$, where $L$ 
is a hyperplane in one of the fibres of $P_Y$. 

To justify the claim, let $F \subset P_Y$ 
be the fibre over the intersection point 
$\{x_0\}=Z \cap Y \subset Y$, let $L$ be the hyperplane 
cut out by $\Ztilde$ on $F$ (actually on $P_Y$), and let 
$\epsilon : X'' \to \Bl_Y X$ be the blowing-up along 
$\Ztilde$. The fibre of $\Ptilde_Y$ above $x$ 
has two components: $\Ftilde$ and $\Xi$. 
Moreover, $\Ptilde_Y$ and $P_{\Ztilde}$ intersect along $\Xi$.  
We first notice that there is a morphism $u : X'' \to X'$.  
Indeed, $(\sigma \circ \epsilon)^{-1}(Z)$ is a Cartier divisor, 
and from the universal property of the blowing-up $\rho$, 
we obtain $u$. Further, the universal property 
is used for $\sigma'$ to imply that the natural  
birational map from $X''$ to $\Bl_{\Ytilde} X'$ 
is defined at any point of $\Xi$ not lying on $\Ftilde$. 
\[
\UseTips
 \newdir{ >}{!/-5pt/\dir{>}}
\xymatrix{ 
 X'' \ar@{..>}[r]^<<<<<<v \ar[dr]^u \ar[d]_\epsilon  &  
     \Bl_{\Ytilde}X' \ar[d]^{\sigma'} \\ 
 \Bl_Y X \ar[d]_\sigma  & X' \ar[d]^{\rho}     \\ 
 X \ar@{=}[r]           & X} 
\] 
We restrict $v$ to $\Ptilde_Y$. Since $\Ftilde$ is  
now a divisor, it follows that the restriction (for which we 
use the same symbol $v$) is defined at the 
generic point of $\Ftilde$, and establishes an isomorphism   
$\Ptilde_Y-\Ftilde \to P_{\Ytilde}-\{x'_0\}$. 
Here $x'_0$ is the point of intersection of $P_{\Ytilde}$ 
with the strict transform of the fibre of $P_Z$ over $x_0$. 
Using the Zariski Main Theorem, we conclude 
that $v$ is a morphism that contracts $\Ftilde$ to $x'_0$. 
The definition of the elementary transforms gives the claim. 

Now by Proposition \ref{2.1}, the vector bundle 
that corresponds to $P_{\Ytilde}$ is determined up to tensoring 
with a line bundle. To finish the proof of the lemma, we apply 
the adjunction formula to obtain that 
$\deg N_{Y|X} = \deg N_{\Ytilde|X'} - 1$.
\qed

\begin{thm} \label{2.3}
Let $X$ be a rationally connected variety and $Y \subset X$ a 
smooth rational curve with ample normal bundle. Then there exists 
a sequence of blowing-ups with smooth $2$-codi\-men\-sional 
centers $\Xtilde \to X$ such that the strict transform 
$\Ytilde$ becomes a quasi-line.
\end{thm}
\proof
We blow-up different well-chosen $2$-codimensional smooth 
subvarieties such that, by Lemma \ref{2.2}, 
the strict transform of $Y$ becomes a quasi-line.
\qed

We end this section by presenting an application of the above theorem.
Let us recall from \cite{BBI} the following definition. A quasi-line 
$Y \subset X$ is called an {\it almost-line} if there is a divisor
$D \in \Div(X)$ such that $D \cdot Y=1$. The main reason for introducing 
this notion was the theorem below, proved in \cite{BBI}:

\begin{thm}[2.1 in \cite{BBI}] \label{2.1BBI} 
Let $X$ be a projective manifold of dimension at least two, $Y \subset X$
a closed, smooth, connected curve with ample normal bundle and 
$Y(1)$ the first infinitesimal neighbourhood of $Y$ in $X$. The following
conditions are equivalent:
\begin{ilist}
  \item The natural restriction map $\Pic(X) \to \Pic(Y(1))$ is
  surjective. 
  \item $Y$ is an almost-line.
\end{ilist}
\end{thm}

For a discussion of the history and motivation of condition (i)
reference \cite{BBI} may be consulted. Theorem \ref{2.1BBI} is
completed by the next proposition, a consequence of Theorem
\ref{2.3}. The proposition shows that, at least birationally, the
situation described in (i) occurs precisely when $X$ is rationally
connected. 

\begin{pro} \label{existenceAl}
Let $X$ be a rationally connected projective manifold. Then there is a
morphism $\sigma :X' \to X$ which is a composition of blowing-ups with
smooth $2$-codimensional centers, such that $X'$ contains an almost-line.
\end{pro}
\proof
We proceed as in the proof of Theorem \ref{2.3}. If $Y \subset X$ is a
smooth rational curve with $N_{Y|X}= \bigoplus_1^{n-1} \Oo_Y(a_j)$, 
$a_1 \leq \cdots \leq a_{n-1}$ and $a_{n-1} > 1$, the quasi-line $Y~$
constructed in loc. cit. is actually an almost-line. Indeed, if $E$ is
the exceptional locus of the last blowing-up, we have $E\cdot Y=1$. If
$a_{n-1}=1$ (\ie $Y$ is already a quasi-line), we take $f:\PP^1\to Y$
to be a degree-b covering, with $b \geq 2$. Applying Theorem II 3.14
in \cite{Ko}, we get that a general deformation of $f$ is an embedding
with ample normal bundle, of numerical type $b_1 \leq \cdots\leq
b_{n-1}$ and such that $b_{n-1} \geq 2$. So the previous argument
applies to give the desired conclusion. 
\qed

The smooth cubic threefold $X\subset \PP^4$ is an example of a
rationally connected manifold that does not contain almost-lines. 
To see this, let $Y\subset X$ be a quasi-line. By the
adjunction formula $-K_X\cdot Y=4$, and since $-K_X\sim 2H$, $Y$ is a
conic. But the Picard group is generated by the hyperplane section 
$H$, hence $Y$ is not an almost-line. After blowing-up one line
$l\subset X$, $\Bl_lX$ becomes a conic bundle, $\pi:\Bl_lX\to
\PP^2$. The pull-back of a general line in $\PP^2$ is a surface $S$
isomorphic to $\PP^2$ blown-up at $6$ points. The surface $S$ contains
a section $Y$ for $\pi|_S$ with self intersection $1$. $Y$ is an
almost-line on $\Bl_lX$.

\section{Rationality via quasi-lines} \label{s3}

In this section we wish to discuss some conditions under which a 
rationally connected manifold is actually rational. 

We start by reviewing a construction that involves the Hilbert scheme
associated to a rational curve on $X$, and the universal family of
this Hilbert scheme; see \cite{Ko}, or \cite{BBI} Section 3. Let $X$
be a projective manifold and $Y \subset X$ a quasi-line. Consider the
Hilbert scheme which corresponds to the Hilbert polynomial (for a
certain polarisation) of $Y$ in $X$. Since  
$H^1(Y,N_{Y|X}) = H^1(\PP^1, \bigoplus_1^{n-1}\Oo_{\PP^1}(1)) = 0$, 
the Hilbert scheme is smooth at $[Y]$, and $[Y]$ lies on a unique
irreducible component, $\Hh$, of this Hilbert scheme. We denote by
$\Yy$ the universal family over $\Hh$, and by $\pi$ and $\Phi$ the two
projections. 
\[
\begin{CD}
  \Yy @>{\Phi}>> X \\
  @V{\pi}VV \\
  \Hh
\end{CD}
\]

For each $h \in \Hh$, the curve $Y_h = \Phi(\pi^{-1}(h))$ satisfies 
$[Y_h]=h$, and the restriction of $\Phi$ to $\pi^{-1}(h)$ is an 
isomorphism onto $Y_h$. In addition, there exists a 
neighbourhood of $[Y]$ in $\Hh$ such that for each closed point 
$h$ in this neighbourhood, the curve $Y_h$ is a quasi-line.

For $x \in Y$, we consider the closed subscheme $\pi (\Phi^{-1}(x))$ 
of $\Hh$; it contains a closed point $h$ if 
and only if the curve $Y_h$ passes through $x$.
If $\Ii_x$ is the ideal sheaf of $x$ in $Y$,
then the tangent space to $\pi (\Phi^{-1}(x))$ at $[Y]$ 
is identified with the space of global sections of 
$N_{Y|X}\otimes \Ii_x$.
Since $N_{Y|X}\otimes \Ii_x \simeq (n-1)\Oo_{\PP^1}$ and $H^1$ 
is trivial, the tangent space is isomorphic to $\CC^{n-1}$ and 
the subscheme is smooth at $[Y]$. 
Let $\Hh_x$ be the unique irreducible component 
that contains $[Y]$. 
We shall use the same notation $\pi$ and $\Phi$ 
for the above projections restricted to the universal family 
$\Yy_x \to \Hh_x$.

\begin{notation*}
  Let $Y$ be a quasi-line on $X$. 
  The number of quasi-lines from the family determined by $Y$ and
  passing through two general points of $X$ will be denoted by
  $e(X,Y)$, cf.~\cite{IoVo}. 
  
  The number of quasi-lines from the family passing through one
  general point of $X$ and tangent to a general tangent vector at that
  point will be denoted by $e_0(X,Y)$. 
\end{notation*}
To see that these numbers are indeed finite, we take $\xi$ a 
$0$-dimensional subscheme of length $2$ in $Y$ in such a way that
$\{x\} \subset \supp (\xi)$. The closed subscheme $\Hh_\xi$ of 
curves through $\xi$ is contained in $\Hh_x$ and, as before, its 
tangent space at $[Y]$ is identified with 
$H^0(Y,N_{Y|X}\otimes \Ii_\xi)$.
This space of global sections is trivial, and this implies the 
finiteness of the number of quasi-lines through $\xi$.
Notice that the degree of $\Phi : \Yy_x \to X$ is equal to $e(X,Y)$.

Recall that a model is a couple $(X,Y)$ with $Y\subset X$ a smooth
rational curve with ample normal bundle.

\begin{pro}  \label{e0e} 
Let $(X,Y)$ be a model with $Y$ a quasi-line. Then the following
assertions hold: 
\begin{ilist}
  \item If $e_0(X,Y)=1$, then $X$ is a unirational variety.
  \item If $e(X,Y)=1$, then $X$ is a rational variety.
\end{ilist}
\end{pro}
\proof
Let $x\in Y$ be a fixed point. Let $\sigma:\Bl_x X\to X$ be the
blow-up of $X$ at $x$ and $E$ be the exceptional divisor. In the
diagram  
\[
\UseTips
  \newdir{ >}{!/-5pt/\dir{>}}
  \xymatrix{   
    && \Bl_x X \ar[d]^\sigma \\
  \Ee\, \arhook{r} & \, \Yy_x \, \ar[r]^\Phi \ar[d]_\pi \ar@{.>}[ur] & \, X \\
  & \Hh_x \ar@/^/[ul]^s}
\]
the map $\Phi$ contracts the divisor $\Ee=\Phi^{-1}(x)$. This
divisor is the image of the natural section $s : \Hh_x \to \Yy_x$ that
maps a point $h \in \Hh_x$ to the point $x$ on the fibre $Y_h$ of the
universal family. Since both $\Hh_x$ and $\pi$ are
generically smooth, the rational map $\sigma^{-1} \circ \Phi$ is
defined at a general point of $\Ee$. It follows that $\sigma^{-1}
\circ \Phi$  maps $\Ee$ to $E$. 

The fact that the restriction of the map $\sigma^{-1}\circ\Phi$ to
$\Ee$ gives a birational isomorphism to $E$ means precisely that
$e_0(X,Y)=1$. As $\Yy_x$ is birationally isomorphic to
$\Ee\times\PP^1$, (i) is proved. 

To see (ii), we first show that the restriction of the rational map
$\sigma^{-1}\circ\Phi:\Ee\dasharrow E$ is dominant. This comes from the fact
that the map is generically finite, since a point 
$y \in \Ee \cap Y$ is sent to the point of $E$ that corresponds 
to the tangent vector of $Y \subset X$ at $x$.
Since $e(X,Y)=1$, $\sigma^{-1}\circ\Phi$ is birational and by the
Zariski Main Theorem, it follows that its restriction to $\Ee$ is also
birational. Thus $\Yy_x$ is rational and so is $X$.
\qed

When $(X,Y)$ is a model with $Y$ a quasi-line, the preceding argument
gives the inequality 
\[
  e_0(X,Y)\leq e(X,Y).
\]
This inequality may be strict. To see how this is so, we use the model
$(X,Y)$, constructed in \cite{BBI} Example 2.7, and some results to be
established in Section \ref{formalGeometry} concerning the behaviour
of $e$ and $e_0$ under \'etale covers along $Y$, see the proof of
Theorem \ref{eEtildeB}. $X$ is the desingularisation of the toric
quotient of $\PP^n$ by the cyclic group of order $n+1$ and $Y$ is a
quasi-line isomorphic to the image of a line in $\PP^n$. 
$e_0(X,Y)$ is preserved by \'etale covers over $Y$ and $e(X,Y)$ is
equal to $e(\PP^n,\text{line})$ multiplied by the degree of the
projection map. Hence $e_0(X,Y)=e_0(\PP^n,\text{line})=1$ and 
$e(X,Y)=n+1$. 

Another example of a model with the distinguished curve a quasi-line for
which the numbers $e$ and $e_0$ can be explicitly computed is given by
the proposition below.

\begin{pro} \label{e0eForConic}
  Let $X\subset\PP^4$ be a smooth cubic threefold and let $\Qq$ be the
  family of all conics lying on $X$. Then:
  \begin{ilist}
  \item $\Qq$ is an irreducible family of dimension $4$.
  \item A conic $\Gamma$ corresponding to a general point of $\Qq$ is
    a quasi-line.
  \item $e_0(X,\Gamma)=e(X,\Gamma)=6$.  
  \end{ilist}
\end{pro}
\proof
It is a classical fact that the family $\Dd$ of lines contained in $X$
is a smooth, irreducible surface and that there are $6$ lines passing
through a general point of $X$, see \eg \cite{Ko}, 266-270. 
Let $G(3,5)$ be the Grassmannian of planes in $\PP^4$. The incidence
$\{(l,P)\mid l\subset P\}\subset \Dd\times G(3,5)$ is irreducible and
has dimension $4$. The projection of this incidence on the second
factor is birational on its image which identifies to $\Qq$. 
Hence (i).

(ii) was first proved by Oxbury in \cite{Ox}. See also \cite{BBI},
Theorem 3.2 for a more conceptual argument.

To prove (iii), we consider $l$ a general line in $\PP^4$ and
$\{x,x',y\}$ the intersection of $l$ with $X$. Clearly, there is a
bijection between lines contained in $X$ and passing through $y$ and
conics contained in $X$ and passing through $x$ and $x'$. Hence
$e(X,\Gamma)=6$. In a similar vein, let $x$ be a general point of $X$,
$v$ a general tangent vector at $x$ and $y$ the other point of
intersection of the line, 
$l_{x,v}$ determined by $v$, with $X$. For every line contained in $X$
and passing through $y$, the plane spanned by this line and $l_{x,v}$
cuts out a residual conic $\Gamma$ through $x$ and tangent to
$v$. This shows that $e_0(X,\Gamma)=6$, too.
\qed

The considerations made at the beginning of this section, together with
Theorem \ref{2.3}, allow us to introduce the definition below.

\begin{definition} \label{eX}
Let $X$ be a rationally connected projective manifold.
We denote by $e(X)$ the minimum of $e(X',Y')$ for all models
$(X',Y')$, where $\sigma :X'\to X$ is a composition of blowing-ups
with smooth centers, and $Y'$ is a quasi-line on $X'$.
\end{definition}

\noindent
The number $e(X)$ leads to a characterization of rational manifolds inside
the class of rationally connected ones.

\begin{thm} \label{birEx}
Let $X$ be a rationally connected projective manifold.
\begin{ilist}
  \item The number $e(X)$ is a birational invariant of $X$.
  \item $X$ is rational if and only if $e(X)=1$.
\end{ilist}
\end{thm}
\proof
(i) Let $\varphi:X_1 \dasharrow X_2$ be a birational isomorphism between
two rationally connected projective manifolds. Let $\sigma :X' \to X_2$
be a composition of blowing-ups and let $Y' \subset X'$ be a quasi-line
such that $e(X_2)=e(X',Y')$. 
Let $\mu = \sigma^{-1}\circ\varphi :X_1 \dasharrow X'$.
By \cite{Hi}, there is $\rho :X \to X_1$ which is a composition of
blowing-ups such that $\mu\circ \rho :X \to X'$ is a birational
morphism. Let $Y \subset X$ be the inverse image by $\mu\circ \rho$ of
a general deformation of $Y'$. We have $e(X_2)=e(X',Y')=e(X,Y) \geq
e(X_1)$. The opposite inequality follows by symmetry.

(ii) Clearly $e(\PP^n)=1$, so if $X$ is rational, then $e(X)=1$ by
(i). The converse follows from Proposition \ref{e0e}.
\qed

The number $e(X)$ seems to be very difficult to compute. On one hand,
its definition involves arbitrary blowing-ups of $X$; on the other
hand, even for fixed $X$, there are in general infinitely many
families of quasi-lines $Y$ on $X$ and it is not at all clear how to
compute the minimum of all numbers $e(X,Y)$. To circumvent this, we
shall introduce in Section \ref{formalGeometry} a ``local version'' of
$e(X)$, denoted $\etilde(X,Y)$. It is associated to a given quasi-line
$Y$ and it is easier to compute. However, it is an open problem
whether or not $\etilde(X,Y)$ leads to a birational invariant. We
refer to the discussion in the last section.

\section{Strongly-rational manifolds} \label{s4}

The aim of this section is to give a convenient characterization of 
strongly-rational manifolds, such that later on we can establish 
their stability with respect to small smooth deformations.
We shall need the following result from \cite{BBI}, 
which is also a particular case of Theorem \ref{thIoVo}. 

\begin{thm}[4.4 in \cite{BBI}] \label{4.1}
$X$ is strongly-rational if and only if $X$ 
contains a quasi-line $Y \subset X$ and a divisor 
$D$ with $D \cdot Y = 1$ and $\dim |D| \geq n$.
\end{thm}
The argument for the converse is by induction on $n =\dim X$. 
It eventually says that the linear system $|D|$ has dimension $n$
and corresponds to the hyperplanes of the projective space.

\begin{thm} \label{4.2}
A manifold $X$ is strongly-rational if and only if

1) $X$ contains a quasi-line $Y$,

2) there exists a point $x \in X$ smooth on every curve 
$Y' \sim Y$ passing through it, and

3) $e(X,Y) = 1$.
\end{thm}
\proof
If $X$ is strongly-rational, then $X$ 
contains an open subset $U$ isomorphic to an open subset $V$
of the projective space, whose complement is at least 
$2$-codimensional. 
Conditions {\it 1)} and {\it 3)} follow by taking $Y$ to be the image 
of a line on $V$.

For {\it 2)}, we consider a general point $x \in U$, a curve $Y$ 
as above, and $\Lambda$ the linear system of divisors that pass 
through $x$ and correspond to hyperplanes of $\PP^n$. Clearly 
$x$ is an isolated base point of $\Lambda$. If this base locus 
reduces to $x$, then, for every $Y \in \Hh_x$, we can choose a divisor 
$D \in \Lambda$ which avoids a point in each irreducible component 
of $Y$. Hence $(D \cdot Y)_x = 1$, for $D \cdot Y = 1$. Thus 
$Y$ is smooth at $x$. To finish the proof, we need to deal 
with the case when $\Bs \Lambda \neq \{x\}$. 
To do this, we consider $\sigma : X' \to X$, 
the blow-up of $X$ along $\Bs \Lambda - \{x\}$ 
and $x' \in X'$ such that $\sigma(x')=x$. A quasi-line 
$Y_0 \subset X$ that does not meet the center of the blowing-up 
gives a quasi-line $Y'_0$ in $X'$, and hence, the Hilbert scheme 
$\Hh'_{x'}$. 
As above, we obtain the result that every $[Y'] \in \Hh'_{x'}$ 
is smooth at $x'$. Now, to descend the result to $X$, the Chow scheme 
will be used; it has the advantage of possesing the
functorial property that comes from the push-forward of cycles, 
unlike the Hilbert scheme.

Let $\Ch_x(X)$ be the irreducible component that corresponds to $Y_0$, 
of the Chow scheme of cycles through $x$, and similarly,
$\Ch_{x'}(X')$ to $Y'_0$. If we take $\widetilde{\Hh}_x \to \Hh_x$ 
to be the normalisation of $\Hh_x$, then the fundamental class of 
an element in $\widetilde{\Hh}_x$ provides us with a natural morphism 
$f: \widetilde{\Hh}_x \to \Ch_{x}(X)$ that is birational and surjective. 
From the diagram
\[
\begin{CD}
  \widetilde{\Hh'}_{x'} @>f'>> \Ch_{x'}(X') \\
  &                  &      @VV{\sigma_\ast}V \\
  \widetilde{\Hh}_x     @>f>>  \Ch_{x}(X) 
\end{CD}
\]
it follows that for any $[Y] \in \Ch_x(X)$, the curve that 
corresponds to it is smooth at $x$. This smoothness and  
Theorem I.6.5 in \cite{Ko} imply that the morphism 
from the universal family to the Chow scheme, $p: \Cc_x \to \Ch_{x}(X)$, 
is smooth at every point of the distinguished section. 
Then, from the commutative diagram,
\[
\begin{CD}
  \widetilde{\Yy}_x @>>>  \Cc_x \\
  @V{\widetilde{\pi}}VV                 @VVpV \\
  \widetilde{\Hh}_x @>f>> \Ch_{x}(X) 
\end{CD}
\]
it follows that $\widetilde{\pi}$ is smooth at every point 
of $\widetilde{\Ee}$, and consequently the same holds for 
$\Yy_x \to \Hh_x$.

\medskip

To prove the converse, we start with some remarks. 
The first one is that the hypotheses give rise to the diagram
\[
\UseTips
 \newdir{ >}{!/-100pt/\dir{>}}
 \xymatrix{
   & \Ee \ar[r]^{\widetilde{\Phi}\mid_\Ee} \ar[d] 
   & E \ar[d]                \\
   & \Yy_x \ar[r]^{\widetilde{\Phi}} \ar[dl]^\pi \ar[dr]_\Phi 
   & \Xtilde \ar[d]^\sigma \\
     \Hh_x \ar@/^/[uur]^s & & X}
\]
with $\sigma$ being the blowing-up of $x$, $E$ the exceptional
divisor, $\widetilde{\Phi}$ a birational morphism and $s$ an
isomorphism. To see that $\widetilde{\Phi}$ is a morphism, it is
sufficient to show that $\Ee$ is a Cartier divisor (see the lemma
hereafter), and to apply the universal property of a blowing-up.

The second remark is that if $F \subset \Yy_x$ is a 
general fibre of $\pi$ and $y \in F -(F \cap \Ee)$ is any point, 
then $\Phi$ is a local isomorphism at $y$.
Indeed, if it is not, from the Zariski Main Theorem, 
$\Phi^{-1}(\Phi(y))$ is positive dimensional. 
Moreover $\Phi(y) = x' \neq x$, hence there are infinitely 
many quasi-lines through $x$ and $x'$, which is impossible. 

The third remark is that there exists an effective divisor
$D \subset X$ such that $D \cdot Y = 1$. Indeed, if 
$H \subset E \simeq \PP^{n-1}$ is a general hyperplane, 
we look at 
$D = \Phi((\widetilde{\Phi}|_\Ee \circ s \circ \pi)^\ast H)$.
Then, the multiplicity of $x$ on $D$ is given by 
\begin{equation} \label{HyperplaneToDiv}
  \mult_x(D) = (-1)^{n-2}(\sigma|E)_\ast(\tilde{D} \cdot E^{n-1}) 
             = H^{n-1} = 1, 
\end{equation}
hence the point $x$ is smooth on $D$. As a consequence of the previous
remark, a general deformation of $Y$ through $x$ meets $D$ only at $x$
and the intersection is transverse. We conclude that $D \cdot Y = 1$.

The fourth remark is that for a general point $p \in E$ and general
hyperplanes $H_1, \ldots , H_{n-1}$ through $p$, if $D_1, \ldots ,
D_{n-1}$ are the divisors on $X$ corresponding to the $H_i$'s as in
(\ref{HyperplaneToDiv}), then the $D_i$'s cut out a quasi-line
transversely. It is sufficient to justify this on $\Yy_x$, since
$\Phi$ is birational. On $\Yy_x$, this is obvious from the generic
smoothness of $\pi$. 

The last remark is that $\dim |D| \geq n$. 
Clearly $\dim |D| \geq n-1$. 
If equality existed, $|D|$ would consist only of divisors 
that come from hyperplanes of $E$. Using the above results, 
the point $x$ would be an isolated base point of the 
linear system $|D|$ and $x$ would be smooth on a general 
divisor. Let $\epsilon : X' \to X$ be the blow-up of $X$ 
along $B$%
\footnote{In fact we sytematically apply Hironaka's result to obtain a 
base-point free linear system through blowing-ups with smooth centers.},
the base locus of $|D|$ minus $x$, and let $|D'|$ be the strict
transform of $|D|$. By Bertini, it would follow that a general divisor
in $|D'|$ is smooth. Let $D_1, \ldots, D_{n-1} \in |D|$ be general
divisors. Their intersection locus would contain a quasi-line $Y$ and
perhaps some other closed subset disjoint from $Y$, and would be 
transverse along $Y$. Now $Y$ is disjoint from $B$, hence 
$Y'=\epsilon^{-1}(Y)$ is a quasi-line on $X'$ and in a neighbourhood of
$Y'$, the strict transforms $D'_1, \ldots, D'_{n-1}$ would also
intersect transversely along $Y'$. From $H^1(X', \Oo_{X'})=0$ and the
exact sequence  
\[
  0 \lra \Oo_{X'} 
    \lra \Oo_{X'}(D') 
    \lra \Oo_{D'_1}(D') 
    \lra 0,
\]
it would follow that $\dim |D'|_{D'_1}| = n-2$, and using 
the exact sequence of normal bundles
\[
  0 \lra N_{Y'|D'_1} \lra N_{Y'|X'} \lra \Oo_{Y'}(1) \lra 0,
\]
that $Y'$ is a quasi-line on $D'_1$. As we have noticed, 
$D'_1$ is smooth, hence $H^1(D'_1, \Oo_{D'_1}) = 0$, and we could 
pursue the restriction procedure. After $n-2$ steps, we 
would arrive at a smooth, regular surface $S$, together with 
a rational curve $Y' \subset S$ such that $(Y')^2 = 1$ and 
$\dim |Y'| \leq 1$. By Riemann-Roch, this would be impossible.

At this point, we invoke Theorem 4.1 and obtain that
$X$ is strongly-rational.
\qed

\begin{lem} \label{4.3}
Let $\Ee \hookrightarrow \Yy \stackrel{\pi}{\to} \Hh$ 
be a flat morphism of schemes of relative dimension $1$, 
with $\Hh$ integral and $\Ee$ a section. 
Then $\Ee$ is a Cartier divisor on $\Yy$ if and only if 
$\Ee$ intersects each fibre in a smooth point of the fibre.
\end{lem}
\proof
The ``only if'' part is clear.
For the converse, let us suppose that for every $e \in \Ee$, the 
local ring $\Oo_{\pi^{-1}(\pi(e)),e}$ is a discrete valuation ring
and let us note $A = \Oo_{\Hh, \pi(e)}$, $B=\Oo_{\Yy,e}$ the local 
rings with $\mmm$ and $\nnn$ the maximal ideals. We have the 
commutative diagram
\[
\UseTips
 \newdir{ >}{!/-5pt/\dir{>}}
 \xymatrix{
    A \ar[r] \ar[rd]^{\simeq} & B \ar[d] \\
                              & B/I}
\]
with $I = \Ii_{\Ee,e}$. Hence, the exact sequence of $A$-modules
\[
  0 \lra I \lra B \lra B/I \lra 0
\]
splits. Tensoring with $A/\mmm$ over $A$, the sequence remains 
exact and gives the injection 
\[
  0 \lra I \otimes_A A/\mmm = I/\mmm I \lra B/\mmm B.
\]
But $B/\mmm B$ is a discrete valuation ring, hence $I/\mmm I$ 
is a principal ideal. The inclusion $\mmm I \subset \nnn I$ 
implies that $I/\mmm I$ surjects onto $I/\nnn I$, that $I/\nnn I$ is 
principal, and applying Nakayama's Lemma, that $I$ is principal.
It remains to be shown that $I$ is generated by a nonzero 
divisor in $B$. Since $B/I \simeq A$ is a domain, 
$I$ is a prime ideal, and the next lemma 
ends the proof.
\qed

We are grateful to Lucian B\u{a}descu for pointing out the next lemma
to the first named author.

\begin{lem} \label{4.4}
Let $B$ be a local Noetherian ring. If $I \subset B$ is a 
principal prime ideal of height $1$, then $I$ is generated 
by a nonzero divisor.
\end{lem}
\proof
Let $I = (b)$ and let $\beta \in B$ such that $\beta b = 0$.
Then either $\beta \notin I$, or $\beta \in I$. 
In the former case, considering the localisation 
$B_I$, we get $b=0$, a contradiction. 
In the second case, $\beta = \beta_1 b$. The same argument 
applied to $\beta_1$ yields $\beta_1 = \beta_2 b$, \ie
$\beta \in I^2$. We may continue the process as long as 
$\beta_n \notin I$ and, if this happens, infer that 
$\beta \in I^n$, for every positive integer $n$.
By the Krull Intersection Theorem, $\beta = 0$.
\qed

\medskip

In the sequel, we shall refer to the set-up 
\[
\UseTips
 \newdir{ >}{!/-5pt/\dir{>}}
\xymatrix{
 X_0   \arhook{r} \ar[d] & \Xx \ar[d]^p \\
 \{0\} \arhook{r}        & B}
\] 
where $\Xx$ is a smooth variety, $B$ a smooth affine
curve and $p$ is a proper and smooth morphism,
as to a small deformation of the smooth projective 
variety $X_0$.

\begin{thm} \label{4.5}
Strongly-rational varieties are stable with respect to 
small deformations.
\end{thm}

Note that such a property is not expected to hold for rational
manifolds. Indeed, several examples of smooth rational cubic four-folds
in $\PP^{5}$ are known to exist, but the general one is expected
to be non-rational.

\proof
We shall show that the three conditions of Theorem \ref{4.2} 
are stable with respect to small deformations. For the first 
one the result is well known, see \eg \cite{BBI} Proposition 3.10, 
and the following lemmas deal with the remaining two conditions. 
\qed

\begin{lem} \label{4.6}
Condition 2) is stable with respect to small deformations.
\end{lem}
\proof
Let $p:\Xx\to B$ be a small deformation of $X_0$. 
$X_0$ contains a point $x$ and a quasi-line $Y$ through 
$x$, such that every other $Y' \in \Hh_x$ is smooth at $x$.

After a base change, if necessary, we can assume that 
$\Xx \to B$ has a section $s$ through $x$. If $\bHh_s$  
is the relative Hilbert scheme over $B$, of quasi-lines through 
some point of $s(B)$, we have the diagram
\[
\UseTips
 \newdir{ >}{!/-5pt/\dir{>}}
 \xymatrix{
   \bEe \ar[r] & \bYy_s \ar[r]^\Phi \ar[d]^\pi & \Xx \ar[dd]^p \\
               & \bHh_s \ar[d]^q \ar@/^/[ul]   &               \\
               & B \ar@{=}[r]                  & B}
\]
which over $0 \in B$, becomes 
\[
\UseTips
 \newdir{ >}{!/-100pt/\dir{>}}
 \xymatrix{
   \Ee_0 \ar[r] & \Yy_0 \ar[r]^{\Phi_0} \ar[d]^{\pi_0} & X_0 \\
                & \Hh_0 \ar@/^/[lu]                        &}
\]
Every fibre of $\pi_0$ is smooth at the point of intersection 
with $\Ee_0$.

Let $\Uu \subset \bYy_s$ be the subset 
of those $y \in \bYy_s$ at which $\pi$ is smooth. Since $\pi$ is flat, 
$\Uu$ is an open subset. The hypothesis yields 
$\Ee_0 \subset \Uu$, and since $q \circ \pi$ is proper, 
$\bYy_s - \Uu$ is sent onto a closed subset of $B$. 
We conclude that $\Ee_b$ is contained in $\Uu$ for every $b$ 
in a neighbourhood of $0$.
\qed

\begin{lem} \label{4.7}
Condition 3) is stable with respect to small 
deformations.
\end{lem}
\proof
We consider a small deformation and the diagram
above. As before, when restricted to $0$, it becomes 
\[
  \begin{CD}
    \Yy_0 @>{\Phi_0}>> X_0 \\
    @VVV                   \\
  \Hh_0
\end{CD}
\]
with $\Phi_0$ a birational morphism. To see that 
$\Phi_b$ is birational in a neighbourhood of $0$, 
one may consider a smooth curve $C \subset \Xx$ through $x$, 
the intersection of general very ample divisors that contain $x$,
and its pre-image $C' \subset \bYy_s$.
The restrictions $p|_C : C \to B$ and 
$(q \circ \pi)|_{C'} : C' \to B$ are proper and quasi-finite,
and hence finite. Consequently $\Phi |_{C'} : C' \to C$ is finite.
But the fibre above $x$ has length $1$ and the result is 
established.
\qed

\section{Quasi-lines and formal geometry} \label{formalGeometry}
 
If $Y$ is a closed subscheme of $X$, the theory of {\it formal
  functions of $X$ along $Y$} was developped by Zariski and
Grothendieck as an algebraic substitute for a complex tubular
neighbourhood of $Y$ in $X$. $\widehat X|_Y$ denotes  the formal
completion of $X$ along $Y$, which is the ringed space with
topological space $Y$ and sheaf of rings
$\Oo_{\Xhat|_Y}=\displaystyle \lim_{\longleftarrow} \Oo_X/\Ii^n$.
$\Ii$ is the sheaf of ideals defining $Y$ in $X$.

In \cite{HM}, Hironaka and Matsumura have introduced and studied
$K(\Xhat|_Y)$, the ring of {\it formal rational functions of $X$ along
  $Y$}. In good cases it is a field that contains the field 
$K(X)$ of rational functions of $X$. We recall the following
definitions from \cite{HM}: $Y$ is G2 in $X$ if $K(\widehat
X|_Y)$ is a field and the field extension 
$K(X)\subset K(\widehat X|_Y)$ is finite. $Y$ is G3 in $X$
if the inclusion $K(X) \subset K(\widehat X|_Y)$ is an isomorphism. 

\begin{notation*}
If $Y$ is G2 in $X$, we denote by $b(X,Y)$ the degree of the field
extension $K(X)\subset K(\Xhat|_Y)$.
\end{notation*}

We also recall the following two results that will be repeatedly used
in the sequel. They are stated in the particular case of quasi-lines,
which is enough for our purposes. Since the normal bundle of
a quasi-line $Y\subset X$ is ample, $Y$ is G2 in $X$, see \cite{H}.

\begin{likethm}[The Hartshorne-Gieseker construction ({\rm see \cite{G},
  Theorem 4.3})]
If $Y\subset X$ is a quasi-line, then there is a morphism of models
$(X',Y')\to(X,Y)$ of degree $b(X,Y)$, \'etale along $Y'$, and such that
$Y'$ is {\rm G3} in $X'$.
\end{likethm}

\begin{likethm}[Gieseker's Theorem ({\rm see \cite{G}})]
Let $(X_i,Y_i)$ be two models with $Y_i\subset X_i$ a quasi-line,
$i=1,2$. Assume that $Y_i$ is {\rm G3} in $X_i$ and
$\Xhat_1|_{Y_1}\simeq\Xhat_2|_{Y_2}$ as formal schemes. Then there are
Zariski open subsets $U_i\subset X_i$ containing $Y_i$ and an
isomorphism $U_1\simeq U_2$ that sends $Y_1$ to $Y_2$.
\end{likethm}

The following definition is similar to Definition \ref{eX}. Here,
for a given quasi-line, we consider \'etale neighbourhoods instead of its
Zariski neighbourhoods.

\begin{definition*}
  Let $Y\subset X$ be a quasi-line. The number $\etilde(X,Y)$ is the
  minimum of $e(X',Y')$, where $X'$ is a projective manifold,
  $Y'\subset X'$ is a quasi-line, $f:X'\to X$ is a generically finite
  morphism, \'etale along $Y'$, and $f(Y')=Y$.
\end{definition*}

The following result shows a useful relationship between the geometry
of quasi-lines and formal geometry. It will play a key role in the sequel.

\begin{thm} \label{eEtildeB}
  If $Y\subset X$ is a quasi-line, then 
\[
  e(X,Y)=\etilde(X,Y)\cdot b(X,Y).
\]
\end{thm}
\proof
We first show that if $f:(X',Y')\to (X,Y)$ is a generically finite
morphism of models, with $Y$ and $Y'$ quasi-lines, \'etale along $Y'$,
then $e(X,Y)=\deg f \cdot e(X',Y')$. Let $y'\in Y'$ be a fixed point
and let $y=f(y')$. For $x\in X$ a general point, we denote by
$x'_1,\ldots,x'_d$ the points of the fibre over $x$, where $d=\deg f$.
We consider the quasi-lines on $X'$ equivalent to $Y'$, that pass
through $y'$ and $x'_i$, for some $1\leq i\leq d$. Their images on $X$
are quasi-lines through $y$ and $x$, and are equivalent to $Y$. The
induced map on Chow schemes $f_\ast: \Ch_{y'}(X') \to \Ch_y(X)$ is
injective when restricted to the open sets parametrising
quasi-lines. This comes from the fact that the considered quasi-lines
on $X'$ do not intersect the ramification divisor of $f$. It follows
that this restriction of $f_\ast$ is also surjective, giving the equality.

Next we choose an $f$ as above, such that
$e(X',Y')=\etilde(X,Y)$. We claim that $Y'$ is G3 in $X'$. We
can apply the Hartshorne-Gieseker construction to get
$g:(X'',Y'') \to (X',Y')$ as above, with $\deg g=b(X',Y')$. Then, by
the previous step,
\[
  \etilde(X,Y)= e(X',Y')= b(X',Y') \cdot e(X'',Y''),
\]
and from the definition of $\etilde(X,Y)$ it follows that
$b(X',Y')=1$.

To finish the proof we consider the following diagram associated to
$f$,
\[
\UseTips
  \newdir{ >}{!/-100pt/\dir{>}}
  \xymatrix{
    K(X) \arhook{r} \arhook{d} & 
         K(\Xhat|_Y) \arhook{d}^\simeq \\
    K(X') \arhook{r}^<<<<<\simeq & K(\widehat{X'}|_{Y'})
  }
\] 
and conclude that $\deg f=b(X,Y)$. Note that the right vertical
isomorphism comes from the fact that $f$ being \'etale along $Y'$,
induces an isomorphism between $\Xhat|_Y$ and $\widehat{X'}|_{Y'}$, see
\cite{G}, Lemma 4.5.
\qed

The next corollary shows that $\etilde(X,Y)$ depends only on the
formal completion of $X$ along $Y$.

\begin{cor}
  Let $(X,Y)$ and $(X',Y')$ be two models with $Y$ and $Y'$
  quasi-lines. If $\Xhat|_Y\simeq\widehat{X'}|_{Y'}$ as formal
  schemes, then $\etilde(X,Y)=\etilde(X',Y')$.
\end{cor}
\proof
We use Hartshorne-Gieseker's construction and suppose that $Y$ and
$Y'$ are G3. As the formal completions are isomorphic, it
follows that $Y$ and $Y'$ have isomorphic Zariski neighbourhoods, by
Gieseker's Theorem. Hence $e(X,Y)=e(X',Y')$. 
\qed

\begin{cor} \label{inequalities}
  Let $Y\subset X$ be a quasi-line. Then
\[
  e_0(X,Y) \leq \etilde(X,Y) \leq e(X,Y).
\]
\end{cor}
\proof
The second inequality comes directly from the theorem. As for the
first one, it is enough to notice, as in the proof of the above
theorem, that $e_0(X,Y)$ is preserved by \'etale covers over $Y$.
\qed

The connection between $e(X,Y)$ and $\etilde(X,Y)$ gives the following
characterisation of the G3 property. 

\begin{cor} \label{etildeG3}
  Let $Y\subset X$ be a quasi-line. $Y$ is {\rm G3} in $X$ if and only if
  $\etilde(X,Y)=e(X,Y)$. In particular, if $e_0(X,Y)=e(X,Y)$, then $Y$
  is {\rm G3}. 
\end{cor}
As a very special case, a quasi-line $Y \subset X$ with $e(X,Y)=1$ is
G3. This generalizes the fact, first noticed by Hironaka, that a line
in the projective space is G3.

It is easy to see that when $X$ has dimension $2$ and $Y\subset X$
is a quasi-line, the formal completion $\Xhat|_Y$ is isomorphic to
$\widehat{\PP^2}|_{\text{line}}$. However, in higher dimensions, the
situation is completely different, as shown by the following example.
See also \cite{IoVo}.
 
\begin{exa}
Let $E$ be a vector bundle over $\PP^2$ associated to the exact sequence
\[
  0 \lra \Oo_{\PP^2}(1) \lra E \lra \Jj_{p,q}(2) \lra 0,
\]
where $p$ and $q$ are two distinct points in $\PP^2$. Note that,
$p$ and $q$ fixed,  $E$ lives in a $2$-dimensional family. 
Let $Y\subset P(E)$ be the quasi-line, lying over a line in $\PP^2$,
given by the construction in \cite{IoVo}, Proposition 4.2. We claim
that for two models of this type, $(P(E),Y)$ and $(P(E'),Y')$, 
the isomorphism $\widehat{P(E)}|_Y\simeq\widehat{P(E')}|_{Y'}$ holds
only if the bundles $E$ and $E'$ are isomorphic.

To justify the claim, we note that by the following lemma,
$e(P(E),Y)=e(P(E'),Y')=1$. Hence, by Corollary
\ref{etildeG3}, $Y$ and $Y'$ are G3 in $P(E)$ and $P(E')$,
respectively. By Gieseker's result, there exist open subsets
$U\subset P(E)$ and $U'\subset P(E')$ and an isomorphism $\phi:U\to U'$,
such that $\phi(Y)=Y'$. Moreover, the complements of $U$ and $U'$ are
at least $2$-codimensional by Lemma 4.4 in \cite{IoVo}, and $\phi$
induces an isomorphism on the Picard groups. The formula for the
canonical class shows that $\phi^\ast\Oo_{P(E')}(1)=\Oo_{P(E)}(1)$. 
Using the exact sequence that defines the vector bundles, we infer
that $\phi$ is an isomorphism outside the fibers over the fixed points
$p$ and $q$. Now, if $D$ and $D'$ are the pull-backs of a line in
$\PP^2$ on $P(E)$ and $P(E')$ respectively, standard computations with
intersection numbers show that 
$\phi^\ast\Oo_{P(E')}(D')=\Oo_{P(E)}(D)$. Here we use that a general
$D$ avoids the indeterminacy locus of $\phi$. It follows that $\phi$
is an isomorphism between the two projective bundles over $\PP^2$. As
the Chern classes of the two vector bundles are the same, the claim is
proved. 
\end{exa}

\begin{lem} 
Let $(X,Y)$ be a model with $Y$ a quasi-line, $E$ be a rank $r\geq 2$
vector bundle over $X$ and $\pi:P(E)\to X$ the canonical projection. If
there exists $Y'\subset P(E)$ a quasi-line that projects
isomorphically onto $Y$, then $e(P(E),Y')=e(X,Y)$. 
\end{lem}
\proof
Let $x'\in Y' \subset P(E)$ and $x=\pi(x')$. Consider the induced
map $\pi_\ast:\Ch_{x'}(P(E))\to\Ch_x(X)$. The open subset of
$\pi_\ast\inv([Y])$ that corresponds to quasi-lines through $x'$ is 
irreducible of dimension $r-1$. This follows from Proposition 4.1. in
\cite{IoVo}: these quasi-lines correspond to lines in a certain
$\PP^r$, passing through a given point. In particular, there is only
one such curve passing through a second point.

Now, let $x''$ be a general point of $P(E)$. The choice of $x''$
implies that any quasi-line
from the family determined by $Y'$ that passes through $x'$ and $x''$
is mapped by $\pi$ to a quasi-line equivalent to $Y$ and passing
through $x$ and $\pi(x'')$. 
\qed

The above lemma allows us to give a similar example with a rank $r$
vector bundle over $\PP^r$. 

\begin{exa} \label{ar}
Let $r \geq 2$, $X=P(T_{\PP^r}^\ast)$ and $Y \subset X$ be an almost-line 
as in \cite{IoVo}, Proposition 4.2. Let also $X'$ be the projective
space of dimension $2r-1$, and $Y' \subset X'$ be a line.
$Y$ is G3 in $X$, but the formal completions $\Xhat|_Y$ and 
$\widehat{X'}|_{Y'}$ are not isomorphic%
\footnote{For $r=2$, this was proved in \cite{IoVo} by an ad-hoc
  argument. Herbert Kurke informed the first named author that he
  independently proved the assertion about the formal completions by a
  completely different method.}. 

Indeed, by the above lemma, $e(X,Y)=1$. Hence $Y$ is G3 in $X$ by
Corollary \ref{etildeG3}. If the formal completions were isomorphic,
by Gieseker's result, $Y$ and $Y'$ would have isomorphic
Zariski neighbourhoods. The complements of these neighbourhoods would be at
least $2$-codimensional by Lemma 4.4 in \cite{IoVo}, hence $X$ and
$X'$ would have isomorphic Picard groups. This is absurd.
\end{exa}

This example is relevant in connection with the following proposition:

\begin{pro} \label{5.3}
Let $X$ be a projective manifold. The following conditions are equivalent:
\begin{ilist}
  \item $X$ is strongly-rational;
  \item $X$ contains a curve $Y$ in such a way that

\indent\indent 
{\rm (a)} $\widehat X|_Y$ is isomorphic (as formal schemes) 
to $\widehat{\PP^n}|_{\text{line}}$, and

\indent\indent 
{\rm (b)} $Y$ is {\rm G3} in $X$.
\end{ilist}
\end{pro}

The equivalence follows from Gieseker's result, using the fact that a
line in $\PP^n$ is G3. 
The Example \ref{ar} shows that at least for $n$ odd, $n\geq 3$, there
are models $(X,Y)$ with $\dim X=n$ that satisfy condition (b), but not
condition (a). On the other hand, the examples constructed in
\cite{BBI}, 2.7 satisfy condition (a), but not (b). So the two
conditions, both very hard to verify in practice, are independent.

\medskip

Our last corollary is a kind of ``formal Torelli theorem'' for
cubic threefolds.

\begin{cor} \label{conicOnCubic}
  Let $X$ and $X'$ be smooth cubic threefolds in $\PP^4$, and let
  $\Gamma\subset X$ and $\Gamma'\subset X'$ be general conics. If
  $\Xhat|_\Gamma\simeq\widehat{X'}|_{\Gamma'}$, then there exists an
  isomorphism $\phi:X\to X'$ such that $\phi(\Gamma)=\Gamma'$.
\end{cor}
\proof
We have noticed in Proposition \ref{e0eForConic} that such a conic is a
quasi-line with $e_0(X,\Gamma)=e(X,\Gamma)=6$. We apply Corollary 
\ref{etildeG3} to deduce that
$\Gamma$ and $\Gamma'$ are G3. By Gieseker's 
result, $\Gamma$ and $\Gamma'$ have isomorphic Zariski
neighbourhoods. As $X$ and $X'$ are Fano manifolds with $b_2=1$, the
birational isomorphism extends to an isomorphism, see \eg \cite{IoVo},
Proposition 1.17.
\qed

\section{Some questions}

A first question arises from the fact that we do not know of any
example of a quasi-line $Y\subset X$ for which
$e_0(X,Y)<\etilde(X,Y)$, cf.~Corollary \ref{inequalities}. Note that
if equality holds, then $b(X,Y)$ is constant in the family of
quasi-lines determined by $Y$. In particular, if $Y$ is G3, then every
$Y'\sim Y$ is G3.

More importantly, we would like to ask the following:

\begin{likethm}[Question]
If $Y\subset X$ is a quasi-line, is $\etilde(X,Y)$ determined by the
field of formal rational functions $K(\Xhat|_Y)$? 
\end{likethm}
More precisely, if
two models have $K(\Xhat|_Y)\simeq K(\widehat{X'}|_{Y'})$ as
$\CC$-extensions, does it follow that $\etilde(X,Y)=\etilde(X',Y')$?
The question may be reformulated as follows: 
If $(X,Y)$ and $(X',Y')$ are two models with $Y$ and $Y'$ quasi-lines
which are G3 in $X$ and $X'$ respectively, and $X$ is
birational to $X'$, is it true that
$e(X,Y)=e(X',Y')$? The equivalence of the two formulations comes from the
Hartshorne-Gieseker construction. We note that a positive answer to
this question would be analogous to the birational invariance of
$e(X)$ established in Theorem \ref{birEx}. However, $\etilde(X,Y)$ is
much easier to compute than $e(X)$.  

A relevant particular case of the above problem concerns rational
manifolds. 

\begin{likethm}[Question]
Let $X$ be a rational projective manifold and $Y\subset X$ be a
quasi-line, {\rm G3} in $X$. Is it true that $e(X,Y)=1$?
\end{likethm}

This question appears as a natural converse of the facts proved in
Proposition \ref{e0e} and in Corollary \ref{etildeG3}: If $e(X,Y)=1$,
then $X$ is rational and $Y$ is G3 in $X$.

Recall from the proof of Corollary \ref{conicOnCubic} that a general
conic $\Gamma$ lying on a smooth cubic threefold $X\subset\PP^4$ is a
quasi-line, is G3 in $X$ and has $e(X,\Gamma)=6$. It follows
that a positive answer to the above question would yield a completely
new proof of the non-rationality of $X$.

The simple example below points
out the difficulty in constructing a counterexample to the above
question. Consider $\Gamma$ a conic in $\PP^3$, a fixed point
$p\in\Gamma$, a general line $l$ and a general smooth curve $C$ both
meeting $\Gamma$. Take $X$ to be the blow-up of $\PP^3$ with center
$p$, $l$ and $C$. The proper transform of the conic becomes a
quasi-line $Y$. Remark that $\etilde(X,Y)$ is independent of the
choice of the curve $C$, the formal completion $\Xhat|_Y$ being the
same. When $C$ is a line, we easily find out that $e(X,Y)=1$. However,
when $C$ is irrational, $e(X,Y)$ must be greater than one. This comes
from the proof of Proposition \ref{e0e}: If $e(X,Y)=1$, the
exceptional divisor lying over $C$ should be rational. Hence, in this
case, $Y$ is not G3 in $X$.

\appendix

\section*{Appendix: a toric example}

If $(X,Y)$ is a model, with $Y$ an almost-line, one may ask the 
following question, related to the hypothesis of Theorem \ref{thIoVo} 
(see \cite{IoVo}, Remark 1.14):
Can we find a linear system $|D|$ on $X$ such that $D \cdot Y=1$ and
$\dim |D| \geq 1$? We reconsider a basic example from \cite{BBI}, 2.7,
in order to prove, via a toric calculation, that the answer to the
above question is, in general, no. 

We refer to \cite{F} for basic notions on toric varieties
and recall here the following facts we shall need about Cartier
divisors on the toric variety $X(\Sigma)$, where $\Sigma$ is
a fan in the lattice $N \subset \RR^n$:

1) The closure of the orbit corresponding to a ray in the fan is an 
irreducible, effective, toric Weil divisor on the variety.

2) Let $D=\sum_i a_iD_i$ be a toric Weil divisor, where $i$ runs  
over the rays in the fan, $D_i$ is the closure of  an orbit 
corresponding to the $i^{th}$ ray, 
$\rho_i \cap N = v_i \ZZ_{+}$, and $a_i \in \ZZ$. 
$D$ is a Cartier divisor if and only if  the map $\psi_D(v_i)=-a_i$ 
can be extended to a piecewise $\ZZ$-linear map on $\Sigma$.

3) The Picard group for a toric variety is spanned by the classes of 
the toric divisors.

4) Let $D$ be a Cartier toric divisor and $\psi_D$ be its associated 
piecewise $\ZZ$-linear map. The dimension of the space of global 
sections of $\Oo_{X(\Sigma)}(D)$ equals the 
number of integer points in the polyhedron 
\[
P_D = \{ u \in N^\ast \otimes \RR \mid u \geq \psi_D \}.
\]

5) Let $f : X(\Sigma_1) \to X(\Sigma_2)$ be a toric map, 
where $\Sigma_i \subset N_i$, $i = 1,2$. 
Giving a toric map is equivalent to giving a lattice homomorphism 
$\sff : N_1 \to N_2$ such that $\sff (\Sigma_1) \subset \Sigma_2$. 
The pull-back of a toric Cartier divisor $D_2$ on $X(\Sigma_2)$ is 
the toric Cartier divisor characterised by the piecewise $\ZZ$-linear map 
$\psi_{D_2} \circ \sff$.

\medskip

The example we are considering is the following:
Let $U_{n+1}$ be the cyclic group of order $n+1$ acting on $\PP^n$ by 
$\zeta \circ [x_0,x_1, \ldots , x_n] = 
[x_0, \zeta x_1, \ldots , \zeta^n x_n]$, 
and let $\pi$ be the quotient map $\PP^n \to X = \PP^n/ U_{n+1}$. 
The projective space and $X$ are toric varieties, 
and $\pi$ is a toric map. To see this, let $N$ be the canonical lattice 
in $\RR^n$ spanned by $e_1, e_2, \ldots , e_n$ over $\ZZ$, and let 
$N' \subset N$ be the sub-lattice spanned by $(n+1)e_1, e_2, \ldots , e_n$. 
The vectors 
$v_1 = (n+1)e_1 - 2e_2 - 3e_3 - \cdots - ne_n$,
$v_i=e_i$ for $i=2,\dots,n$ and
$v_{n+1}=-(n+1)e_1+e_2+2e_3+\cdots+(n-1)e_n$ 
form the fans $\Sigma_{N'} \subset N'$ and $\Sigma_N \subset N$, 
with maximal cones the $n+1$ simplicial cones. 
Then the projective space is the toric variety $X(\Sigma_{N'})$ 
and $X = X(\Sigma_N)$. The quotient $N/N' \simeq U_{n+1}$ 
acts naturally on $X(\Sigma_{N'})$, the action being the one considered 
above. The toric map induced by the inclusion 
$\Sigma_{N'} \subset \Sigma_N$ is the quotient map. 

Let $H \subset \PP^n$ be the hyperplane corresponding to the piecewise 
$\ZZ$-linear map $\psi_H (v_j) = - \delta_{2j}$ defined on $\Sigma_{N'}$. 
This map is not the restriction of a piecewise $\ZZ$-linear map 
on $\Sigma_N$. But, if we consider the open subset $U \subset X$ 
corresponding to the sub-fan 
$\Delta \subset \Sigma_N$, where $\Delta$ is the union of all the rays 
in $\Sigma_N$, then $\psi_H$ defines a divisor $D$ on $U$. 
We note by $\psi_D = \psi_H |_{\Delta}$ the map associated to $D$.

\begin{likethm}[Lemma A.1]
$h^0(U, \Oo_U(D)) = 1$.
\end{likethm}
\proof
We have to count the number of integer points in the polyhedron $P_D$. 
Since 
\[
  P_D \cap N^\ast = 
    \{ u = \sum_{j=1}^{n} \alpha_j e_j^\ast   \mid u \geq \psi_D \} 
    \cap N^\ast
\]
and $u \geq \psi_D$ is equivalent to 
\[
\left\{
  \begin{array}{rcrcl}
     (n+1)\alpha_1 &-& \sum_2^n j\alpha_j     &\geq & 0 \\
                   & &               \alpha_2     &\geq &-1 \\
                   & &               \alpha_j     &\geq & 0 
                           \quad \text{for every $j\geq 3$} \\ 
    -(n+1)\alpha_1 &+& \sum_2^n (j-1)\alpha_j &\geq & 0,
  \end{array}
\right.
\]
we conclude that $P_D \cap N^\ast = \{ 0 \}$.
\qed

\begin{rem*}
The open subset $U \subset X$ equals $X$ minus its $n+1$ singular points 
and $\pi^{-1}(U)$ is the projective space minus the fixed points for the 
$U_{n+1}$-action. 
Moreover 
$\pi^\ast \Oo_U(D) = \Oo_{\PP^n}(H) |_{\pi^{-1}(U)}$, and since 
\[
H^0(\pi^{-1}(U), \Oo_{\PP^n}(H) |_{\pi^{-1}(U)}) = 
H^0(\PP^n, \Oo_{\PP^n}(H)) 
\]
as $U_{n+1}$-representations, the space of global sections for $\Oo_U(D)$ 
is isomorphic to the $U_{n+1}$ trivial sub-representation of 
$H^0(\PP^n, \Oo_{\PP^n}(H))$ which is $1$-dimensional. 
\end{rem*}

Let $\Xtilde \to X$ be a toric desingularization of $X$ obtained by 
taking a smooth sub-division $\widetilde{\Sigma}$ of $\Sigma_N$.

\begin{likethm}[Lemma A.2]
If $\psi_{\Dtilde}$ is a piecewise $\ZZ$-linear extension
of $\psi_D$ to $\widetilde{\Sigma} \subset N$, then the space of global 
sections of $\Oo_{\Xtilde}(\Dtilde)$ is of dimension $\leq 1$.
\end{likethm}
\proof
This is clear, since from $\Sigma_N \subset \widetilde{\Sigma}$ and 
$\psi_{\Dtilde} |_{\Sigma_N} = \psi_D$, we infer that 
$P_{\Dtilde} \subset P_D$. 
\qed

Going back to the problem of finding effective divisors $\Dtilde$ on 
$\Xtilde$ such that $\Dtilde\cdot g^\ast Y=1$, where $Y$ is the 
almost line $\pi(L) \subset U$, 
$L = \{ x_0 = x_1, \, x_2= x_3, \, x_j = 0 \text{ otherwise} \}$, 
the answer is the following:

\begin{likethm}[Corollary A.3]
If $\Dtilde \subset \Xtilde$ is a divisor such that 
$\Dtilde \cdot g^\ast Y=1$, then 
$h^0(\Xtilde, \Oo_{\Xtilde}(\Dtilde) \leq 1$.
\end{likethm} 
\proof
We have the following commutative diagram, 
\[
\begin{CD}
\widetilde{\PP^n} @>f>> \PP^n \\
@V{\widetilde{\pi}}VV @VV{\pi}V \\
\Xtilde @>g>> X
\end{CD}
\]
where $\widetilde{\PP^n}$ is the toric variety 
$X(\widetilde{\Sigma}_{N'})$, with $\widetilde{\Sigma}_{N'}$ the fan 
spanned by $\widetilde{\Sigma}$ in $N'$. Then
\[
n = \widetilde{\pi}^\ast (\Dtilde \cdot g^\ast Y )
  = \widetilde{\pi}^\ast \Dtilde \cdot nf^\ast L ,
\]
hence $\widetilde{\pi}^\ast \Dtilde \cdot f^\ast L=1$. 
Since $f_\ast f^\ast L=L$, it follows that
$f_\ast\widetilde{\pi}^\ast\Dtilde=H$ and  
that the function $\psi_{\widetilde{\pi}^\ast \Dtilde}$ 
defined on $\widetilde{\Sigma}_{N'}$ should send 
to $1$ exactly one of the $n+1$ vectors $v_j$, $j=1,\ldots,n+1$, 
and to $0$ the remaining $n$. Hence, the piecewise $\ZZ$-linear map 
$\psi_{\Dtilde}$, modulo an $SL(n,\ZZ)$ transformation, is 
the map in lemma A.2, and the result follows from that lemma.
\qed

\paragraph{Acknowledgement}
The first named author is grateful to the Department of Mathematics of
the University of Angers for several pleasant and fruitful visits. He
also acknowledges partial financial support from the grant EURROMMAT
of the EU.

\begin{flushleft}
Paltin {\sc Ionescu} \\
Department of Mathematics \\
University of Bucarest \\
14, Academiei Street \\
RO-70109 Bucharest, Romania\\
 {\it and} \\
Institute of Mathematics of the Romanian Academy\\
P.O.Box 1-764\\
RO-70700 Bucharest, Romania \\
{\it e-mail: Paltin.Ionescu@imar.ro}
\end{flushleft}

\begin{flushleft}
Daniel {\sc Naie} \\
Department of Mathematics \\
Universit\'e d'Angers \\
2, bd. Lavoisier \\
FR-49045 Angers, France \\
{\it e-mail: Daniel.Naie@univ-angers.fr}
\end{flushleft}

\end{document}